\documentclass[11pt]{article}
\usepackage{amssymb,amsmath}
\usepackage{graphicx}
\newtheorem{Lemma}{Lemma}
\newtheorem{Theorem} {Theorem}
\newtheorem{Proposition} {Proposition}
\newtheorem{Corollary} {Corollary}

\newtheorem{Definition}{Definition}
\newtheorem{Example}{Example}
\newcommand{\bZ}{{Z}}
\newcommand{\cH}{{\cal H}}
\newcommand{\cK}{{\cal K}}
\newcommand{\haar}{{\hbar}}
\newcommand{\what}{\widehat}
\newcommand{\mbC}{\mathbb{C}}
\newcommand{\mbN}{\mathbb{N}}
\newcommand{\mbR}{\mathbb{R}}
\newcommand{\mbZ}{\mathbb{Z}}
\newcommand{\mcU}{\mathcal{U}}
\newcommand{\mcV}{\mathcal{V}}
\newcommand{\borel}{{\mathcal{B}}}
\newcommand{\mfs}{{\mathbf{s}}}
\newcommand{\mft}{{\mathbf{t}}}
\newcommand{\mfB}{{\mathbf B}}
\newcommand{\mfK}{{\mathbf K}} 
\newcommand{\mfR}{{\mathbf R}}
\newcommand{\defin}{:=}
\newcommand{\ds}{}
\newcommand{\cha}[1]{\left\langle #1 \right\rangle}
\newcommand{\lin}[1]{\overline{\mbox{span}}\left\{ #1 \right\} }
\newcommand{\proof}{\emph{Proof}. }
\newcommand{\qed}{\hfill $\blacksquare$}

\begin{document}


\title{Spectrum of Periodically Correlated Fields}
\author{Dominique Dehay \\
{\small \it Institut de Recherche Math\'ematique de Rennes, CNRS umr 6625, } \\
{\small \it Universit\'e Rennes 2, cs 24307, 35043 Rennes, FRANCE; dominique.dehay@univ-rennes2.fr}  \\
                     Harry Hurd \\
{\small \it Department of Statistics, University of North Carolina,} \\
{\small \it Chapel Hill, NC 27599-2630, USA; hurd@stat.unc.edu } \\
Andrzej Makagon
\thanks{\scriptsize The paper was partially written during the author's stay at Universit\'{e} Rennes 2 - Haute
Bretagne, Rennes, France, in June 2011.}\\
{\small \it Department of Mathematics, Hampton University}\\
{\small \it Hampton, VA 26668, USA; andrzej.makagon@hamptonu.edu}
}
\date{\today}
\maketitle
\begin{abstract}
The paper deals with Hilbert space valued fields over any locally compact Abelian group $G$, in particular over $G=\mbZ^n\times \mbR^m$, which are periodically correlated (PC) with respect to a closed subgroup of $G$. PC fields can be regarded as multi-parameter extensions of PC processes.
We study structure, covariance function, and an analogue of the
spectrum for such fields. As an example a weakly PC field over
$\mbZ^2$ is thoroughly examined.
\end{abstract}


\maketitle

\section{Introduction}
Periodically correlated (PC) processes and sequences have been
studied for almost half of the century and at present they are
very well understood mainly due to works of Gladyshev
\cite{gladyshev61,gladyshev63}, Hurd
\cite{hurd69,hurd74,hurd74b,hurd89b,hurdk91} and other authors
\cite{dehay94,honda82,makms94,mak99,mak00,miamee90,miamees79}. A
summary of the theory of PC sequences can be found in \cite{HM}.
Surprisingly, there are only several works \cite{Bose98,
Bose01,Chen96,dehayhurd,DragJav82,Gardner06,gaspar04,hurdkf,serp05}
dealing with PC fields, and each one concentrates on a particular
type, namely coordinate-wise strong periodicity. An intention of this paper is to  sketch a unified theory of fields over any locally compact Abelian (LCA) group $G$ which are periodically correlated with respect to an arbitrary
closed subgroup $K$ of $G$. We emphasize the case of $G=\mbZ^n\times\mbR^m$
to illustrate the results.
This work includes  stationary fields as well the weakly periodically correlated fields, that is the fields whose covariance function exhibits periodicity (or stationarity) in fewer directions than the dimension of the group.  In the latter case we assume a certain integrability condition (see Definition~\ref{def:square-integr}) in order to develop some simple spectral analysis of those fields.
A work in progress treats the case where this condition is not satisfied.

The paper is organized as follows. In the remaining part of this section we introduce  notation and vocabulary used in the paper, review needed
facts from harmonic analysis on LCA groups, and outline the theory of one-parameter PC processes. In the next three sections we study the covariance function, the notion of the spectrum, and the structure of a $K$-periodically correlated field. These sections include the main results of the paper (Theorems \ref{thm-Corr}, \ref{thm-spectrum} and \ref{PCF-structure}). The last section contains examples that illustrate the theory developed. In particular Example~\ref{ex-wk} gives a complete analysis of the weakly periodically correlated fields over $\mbZ^2$,  introduced in \cite{hurdkf}.

\subsection*{Background}
To avoid confusion and to set the notations of the paper we recall some  features of group theory, Haar measures, Fourier transform, and periodic functions. For more information
on these subjects the authors refer to~\cite{HR,reiter,rudin}.

\bigskip
1.\,\,\underline{Quotient groups, cross-sections, Haar measure, and Fourier transform}\,.\,
Let $G$ be an additive 
locally compact Abelian (LCA) group, $\what{G}$ be its
 dual (group of continuous characters), and let $\cha{\chi,t}$  denote the value of a character $\chi \in \what{G}$ at $t\in G$.   The dual $\what{G}$ can be given a topology that makes it an LCA group such that $\what{(\what{G})} = G$.  Let  $K$ be  a closed subgroup of $G$. The symbol $G/K$ will stand
for the quotient group and $\what{(G/K)}$ for its dual. Let
$\imath$ denote the natural homomorphism of $G$ onto $G/K$, $\imath(t) \defin t+K$, and $\imath^*$ be its dual map $\imath^*:\what{G/K}\to\what{G}$, defined as
$\cha{\imath^*(\eta),t} = \cha{\eta , (t+K)}$ for $\eta\in\what{G/K}$ and $t\in G$. The
mapping $\imath^*$ is injective and continuous,  and for each
$\eta \in \what{G/K}$, $\cha{\imath^*(\eta),\cdot}$ is a
$K$-periodic function on $G$ (see below).
Consequently $\what{G/K}$ can be identified with a closed subgroup $\Lambda_K$
of $\what{G}$ consisting of  the elements $\lambda \in
\what{G}$ such that $\cha{\lambda,t} =1$ for any $t\in K$.
In the sequel we  use the notation $\cha{\lambda,t}\defin\cha{\lambda,\imath(t)}$,
for all $\lambda\in\Lambda_K$ and $t\in G$.
 By $\borel(G)$ we denote the $\sigma$-algebra of Borel sets on $G$.
A \emph{cross-section} $\xi$ for $G/K$ is a mapping $\xi: G/K\to G$ such that

\smallskip\noindent
\begin{tabular}{rll}
\vspace{1.1mm}
&(i)&$\xi$ is Borel,\\ \vspace{1.1mm}
&(ii)&$\xi(G/K)$ is a measurable subset of $G$,\\ 
&(iii)&
$\xi(0) = 0$ and $\xi\circ\imath(t) \in t+K$ for all $t\in G$, where $t+K\defin\{t+k\!:\!k\in K\}$.
\end{tabular}

\smallskip\noindent
For existence and other
properties of a cross-section please see \cite{EK:84,varadarajan}. For each cross-section $\xi$ for $G/K$, the sets
$k+\xi(G/K)$, $k\in K$, are disjoint and their union is $G$, and
hence each element $t\in G$ has a unique representation $t = k(t)
+ \xi(\imath(t))$, where $k(t)\in K$. Note that the function $\xi$ is not
additive, that is $\xi(x+y)$ may be different than $\xi(x)+
\xi(y)$, $x,y \in G/K$.

Any LCA group has a nonnegative translation-invariant measure, unique up to a multiplicative constant, called a Haar measure.
The Haar measures on $G$ and $\what{G}$ can be normalized
in such a way that the following implication  holds
  
\begin{gather*}
\mbox{if}\quad  f\in L^1(G),\quad\what{f}(\chi) \defin
\int_{{G}} \cha{\chi,t} f(t)\,\haar_G (dt)\quad \mbox{for}\,\,\chi \in \what{G},\quad\mbox{and}\quad \what{f}\in L^1(\what{G}) \quad\\
\mbox{then} \quad  f(t) = \int_{\what{G}} \overline{\cha{\chi,t}} \what{f}(\chi)
\,\haar_{\what{G}} (d\chi) \quad \mbox{for a.e.}\,\, t\in {G}. 
\end{gather*}
The function $\what{f}$ above is called the \emph{Fourier transform} of $f$. Here and what follows $L^1(G)$ stands for the space of complex functions on $G$ which are integrable with respect to $\haar_G$, and $\haar_G$   denotes the
normalized Haar measure on the group indicated in the subscript. Note that the normalization of the Haar measures of $ G$ and $\what{G}$ is not unique. We  follow the usual convention that if $G$ is compact and
infinite then the normalization is such that $\haar_G(G) = 1$;
if $G$ is discrete and infinite then the normalized Haar measure of any single
point is 1; if $G$ is both compact and finite then its dual is
also and the Haar measure on $G$ is normalized to have a mass 1 while the
Haar measure on $\what{G}$ is counting measure. The normalized Haar measure on $\mbR$ is the Lebesgue
measure divided by $\sqrt{2\pi}$. Finally, if $K$ is a closed subgroup of $G$ then the normalized Haar measures satisfy Weil's formula
\begin{equation}\label{eq:Weil}
\int_{G/K}\left(\int_K f(k+s)\,\haar_K(dk)\right)\haar_{G/K}(d\dot{s})=\int_G f(t)\,\haar_G(dt), \quad  f\in L^1(G).
\end{equation}
The inner integral above depends only  on the coset $\dot{s}\defin s+K$. See e.g. \cite[Section\,III.3.3]{reiter}.

If  $f\in L^1(G)$ then  $\what{f}$ is a continuous bounded function on $\what{G}$ but not necessarily integrable.
The Fourier transform, which is customarily denoted by the integral $\what{f}(\chi) = \int_{{G}} \cha{\chi,t} f(t)\,\haar_G (dt)$ (even if $f$ is not integrable)
extends from $L^1(G)\cap L^2(G)$ to an isometry
from $L^2(G)$ onto $L^2(\what{G})$  (Plancherel theorem~\cite{rudin}).
If there is a danger of confusion we will recognize the difference by writing
\begin{equation*}
\what{f}(\chi) \stackrel{L^2}{=} \int_{{G}} \cha{\chi,t} f(t)\,\haar_G (dt),\quad  f\in L^2(G).
\end{equation*}

In the sequel we  say that \emph{the inverse formula holds for $f$}   if the function $f$ is the inverse Fourier transform of $\what{f}$. If both $f$ and $\what{f}$ are integrable then clearly the inverse formula holds for both.  Also if $G$ is discrete and $f\in L^2(G)$, then the inverse formula holds for $f$. Indeed, in this case  $\what{f} \in L^1(\what{G})$ because $\what{G}$ is compact, and hence $f(t) = \int_{\what{G}} \overline{\cha{\chi,t}} \what{f}(\chi)\,
\haar_{\what{G}} (d\chi)$  for all $t\in {G}$.

For a separable Hilbert space $\cH$ with inner product
$(\cdot,\cdot )_\cH$ and norm $\|\cdot\|_\cH$, let $L^p(G;\cH)\defin L^p(G,\haar_G;\cH)$, $p=1$ or $2$, be the space of $\cH$-valued fields on $G$ which are $p$-integrable with respect to Haar measure $\haar_G$, that is, $f\in L^p(G;\cH)$ means that $f:G\to\cH$ is $\haar_G$-measurable and the real-valued function $t\mapsto \|f(t)\|_\cH^p$ is integrable with respect to $\haar_G$. It is well known that the space $L^1(G;\cH)$ is a Banach space with the  norm
$\|f\|_{L^1}\defin\int_G\|f(t)\|_{\cH}\,\haar_G(dt),\, f\in L^1(G;\cH)$,
 and the space $L^2(G;\cH)$ is a separable Hilbert space with the inner product
$\big(f,g\big)_{L^2}\defin\int_G                     \big(f(t),g(t)\big)_{\cH}\,\haar_G(dt),\, f,g\in L^2(G;\cH).$
See e.g.~\cite[Chapter III]{DS} (see also~\cite{Dinc2000,Hille48,RudinFA}).
Whenever $f\in L^1(G;\cH)$ then $f$ is Bochner integrable (also called strongly integrable) and its Fourier transform exits. 
Futhermore Plancherel theorem applies and defines an isometry from $L^2(G;\cH)$ onto $L^2(\what{G};\cH)$ (one-to-one), so $\what{f}\in L^2(\what{G};\cH)$ is also well defined for $f\in L^2(G;\cH)$.

\bigskip
2.\,\,\underline{Periodic functions}\,.\,
 Given  $G$ and a closed subgroup
$K$ of $G$, it is natural to call a function $f$ defined on $G$
to be \emph{$K$-periodic} if
\begin{equation*} 
f(t+k) = f(t)\quad\mbox{for all}\,\,t\in G\,\,\mbox{and}\,\,k\in K.
\end{equation*}
In this case, 
the function $f$ is constant on cosets of $K$. Hence a function $f$ on $G$  is
$K$-periodic if and only if  $f$ is of the form $f = f_K\circ \imath $, where
$f_K$ is a function on $G/K$. The
concrete realization $\Lambda_K\defin\imath^*(\what{G/K})\subset\what{G}$ of \  $\what{G/K}$ as a subgroup of $\what{G}$ will be in the sequel called \emph{the domain of the spectrum of $f$}. Note that $\Lambda_K$ is not determined uniquely by $f$, for a
$K$-periodic function can be at the same time periodic with respect to a larger subgroup
$K'\supset K$; in other words we will not be assuming that $K$ is the "smallest" period of $f$.

If  $f\in L^1(G/K;\cH)$, $\cH$ being the set of complex numbers $\mbC$ or any separable Hilbert space,   we consider the Fourier transform of $f_K$ at $\lambda \in
\Lambda_K$
 \begin{equation} \label{f_K-lambda}
  \what{f_K}(\lambda) \defin \int_{{G/K}} \cha{\lambda , x}
f_K(x)\, \haar_{G/K} (dx)\end{equation}
that will be referred to as the \emph{spectral coefficient of  $f$ at frequency} $\lambda \in\Lambda_K$.

A couple of remarks regarding the above definition and its relation to the standard notions of
the spectrum and its domain are certainly due here.  The word {\em spectrum} comes originally from physics,
operator theory, and more recently from signal
processing. It is widely used in the theory of second order stochastic processes. Intuitively,
the spectrum of a scalar function $f$ is a Fourier transform of $f$ in whatever sense it exists. If
$f$ is a locally integrable function on $G=\mbZ^n\times\mbR^m$  
then the spectrum $F$  of $f$ is a Schwartz distribution on $\what
{G}$, which is a functional on a certain space of functions on $\what{G}$ determined by the
relation  $\ds F\big(\what{\phi}\big) = \int_G {\phi}(t) f(t)\, \haar_{G} (dt)$, where $\phi$ runs over the
set of compactly supported functions on $G$ which are infinitely many times differentiable in
last $m$ variables. One can show that if $f$ is additionally  $K$-periodic, then the support
of  $F$ (as defined in \cite{RudinFA}) is a subset of $\Lambda_K$. This rationalizes the name ''\emph{domain of the spectrum}'' that we have
assigned for $\Lambda_K$, as well the phrase ''\emph{the spectrum  sits on}  $\Lambda_K$'' which we will use sometimes.
The first task in understanding the spectrum of a $K$-PC field is thus
to identify the domain of its spectrum or its second order spectrum. (See  below).

The coefficient  $\what{f_K}(\lambda)$ defined in~(\ref{f_K-lambda}) represents an "amplitude" of the
harmonic $\cha{\lambda, \cdot}$ in a spectral decomposition of $f$.
Indeed, if  $f_K$ and $\what{f_K}$ are integrable, then
$\ds f_K(x) = \int_{\Lambda_K} \overline{\cha{\lambda,x}} \, \what{f_K}(\lambda)
\,\haar_{\Lambda_K} (d\lambda)$, $x\in G/K$,
and as a consequence of Weil's formula~(\ref{eq:Weil}) and the fact that
$ \cha{\lambda , t}=\cha{\lambda , \imath(t)} $, $t\in G$, $\lambda \in \Lambda_K$,  we conclude
that
\begin{equation} \label{FTofPF}
f(t) = \int_{\Lambda_K}
\overline{\cha{\lambda, t}} \, \what{f_K}(\lambda)\, \haar_{\Lambda_K}
(d\lambda)
, \quad  t\in G.
\end{equation}
If $\what{f_K}$ is not integrable, then equality~(\ref{FTofPF}) holds only for $\haar_G$-almost every $t\in G$ or is not valid as
stated, but  ${a}_{\lambda}$ still retains its interpretation.

\bigskip
For illustration suppose that $f$ is a continuous  scalar function on $\mbR$ which is periodic
with period $T>0$,  that is such that $f(t) = f(t+T)$ for every $t\in \mbR$. In this case $G=
\mbR$, $K=\{kT\!:\! k\in \mbZ\}$, the
quotient group $G/K$ can be identified with $[0,T)$ with addition
\emph{modulo $T$}, the mapping $\imath$ is defined as \
$\imath(t) = \big[t\big]_{T}$, the remainder in integer division of $t$ by $T$, and the identity $
\xi(x) = x$, $x\in[0,T)$, is the most natural cross-section for $G/K$. The function $f_K$ is
defined as
$f_K(x) = f(\xi(x)) = f(x)$, $x\in [0,T)$. The dual of
$G/K$ is identified with  the subgroup  $\Lambda_K=\{2\pi j/T\!:\! j\in \mbZ\}$ of $\mbR$, and with this
identification $\cha{\lambda , \imath(t)} = \cha{\lambda , t} =  e^{-i\lambda t}$, $\lambda
\in \Lambda_K$, $t \in \mbR$. The normalized Haar measures on $[0,T)$ and $\Lambda_K$ are the
Lebesgue measure divided by $T$ and the counting measure, respectively.  The domain of the
spectrum of $f$ is therefore the set $\Lambda_K$. The spectral
coefficient $\hat{f}_j$ of $f$ at $\lambda = 2\pi j/T$, is given by
$$\ds \hat{f}_{j}=\frac{1}{T}\int_0^T e^{-i2\pi jt/T}f(t)\,dt,\quad j\in\mbZ.$$
Note that the sequence $\{\hat{f}_j\}$ is square-summable and consequently
$\ds  f(t) = \sum_{j= -\infty}^{\infty} e^{i2\pi j t/T} \hat{f}_{j}$, where the series above
converges in $L^2[0,T]$, so in $L^2([-A,A])$ for every $0<A<\infty$.
 The spectrum $F$ of $f$ is  defined by the relation $\ds  F\big(\what{\phi}\big) = \frac{1}{\sqrt{2\pi}}
\int_{-\infty}^{\infty}{\phi}(t) f(t)\,dt
 $. If $\phi$  is an infinitely times differentiable with compact support then
\begin{equation*} F\big(\what{\phi}\big) =  \frac{1}{\sqrt{2\pi}} \int_{-\infty}^{\infty} {\phi}(t) f(t)\,dt 
= \int_{\mbR} \what{\phi}(t) F(dt),
\end{equation*}
where $\ds F = \sum_{j= -\infty}^{\infty} \hat{f}_j\, \delta_{\{2\pi j/T\}}$, and $\delta_a$ denotes the
measure of mass 1 concentrated at $\{a\}$. The spectrum  of $f$ can be therefore identified
with a $\sigma$-additive complex measure $F$ on $\mbR$ sitting on   $\Lambda_K$ and defined by
$\ds F = \sum_{j= -\infty}^{\infty} \hat{f}_j\, \delta_{\{2\pi j/T\}}$.  If the sequence $\{\hat{f}_j\}$ is summable,
then $F$ is a finite measure, but it does not have to be in general.

\section{Periodically Correlated Fields}
Let  $\cH$ be a separable Hilbert space with the inner product
$(\cdot,\cdot )_\cH$. In a probabilistic context the space $\cH$
represents the space of zero-mean complex random variables with
finite variance. 
A \emph{(stochastic) field} $X=\{X(t)\!:\!t\in G\}$
is a measurable function $X: G\to \cH$. 
Let $\cH_X  \defin \lin{X(t)\!:\! t\in G}$  be the smallest closed linear subspace of $\cH$ that contains all $X(t)$, $t\in G$.  The function $\mfK_X(t,s) \defin \big(X(t), X(s)\big)_{\cH}$, $t,s\in G$, is referred to as the \emph{covariance function} of the
field $X$.
A  field $X$ is called \emph{stationary} if it is continuous and for all $t,s\in G$, the function $\mfK_X(t+u,s+u)$
does not depend on $u\in G$. If $X$ is stationary then
$\mfK_X(t+s,s)=\mfK_X(t+0,0) =:\mfR_X(t)$, $t,s \in G$,
and $\mfR_X$ has the form
\begin{equation*}
\mfR_X(t) = \int_{\what{G}}\overline{\cha{\chi,t}}\,\Gamma(d\chi), \end{equation*}
where $\Gamma$ is a non-negative Borel measure on
$\what{G}$~(Bochner Theorem \cite[Section\,IV.4.4]{reiter}). A  field $X$ is called \emph{harmonizable} if
there is a (complex) measure $\digamma$ on $\what{G}\times\what{G}$
such that
\begin{equation} \label{harm}
\mfK_X(t,s) =
\int\!\!\!\int_{\what{G}^2} \overline{\cha{\chi,t}} \cha{\beta,s}
\,\digamma (d\chi, d\beta),\quad t,s\in{G},
\end{equation}
see \cite{rao, rao2005}. 
The measure $\digamma$ above is  called the   \emph{ second order spectral (SO-spectral) measure} of the harmonizable field
$X$. Note that every stationary field is harmonizable with the measure $\digamma$ sitting on
the diagonal: $\digamma (\Delta) =\Gamma\{\chi\in \what{G}\!:\!(\chi,\chi) \in \Delta\}$, $\Delta
\in \borel(\what{G}\times\what{G})$.
The SO-spectrum of a harmonizable field $X$  is the spectral measure $\digamma$ associated with
function $\mfK_X(t, -s)$, $ s,t \in G$ via relation~(\ref{harm}).  Here 
we adopt the terminology of ''\emph{second order spectrum (SO-spectrum)}'' of the field $X$, instead of the usual term ''spectrum'', in order to avoid confusion with the spectrum of a periodic function (or field). By this way we point at the fact that we are considering not the field $X$ by itself but its  covariance function $\mfK_X$. This leads to the following
definition.
\begin{Definition} \label{def:spectrum of X}  Let $X$ be a continuous stochastic field over $G
$. The \emph{second order spectrum (SO-spectrum)} of the field $X$ is the spectrum (Fourier transform, in whatever sense it may exist) of the function
$ G\times G \ni (t,s) \mapsto \mfK_X(t, -s). $
The domain of the SO-spectrum of the field $X$ is defined as the domain of the spectrum of this function.
\end{Definition}
\begin{Definition} 
Let  $K$ be  a closed subgroup of $G$. A field $X$ is called
\emph{$K$-periodically correlated ($K$-PC)} if $X$ is  continuous 
and  the function $G\ni u \mapsto \mfK_X(t+u,s+u)$ is $K$-periodic in $u$ for  all $t,s\in G$.
The group $K$ will be called the period of the PC process $X$.
\end{Definition}
If $G=\mbR$ (or $\mbZ$) and $K=\{kT\!:\! k\in\mbZ\}$ then we will use the phrase ''\emph{PC process
(or  PC sequence) with period $T>0$}\,'', rather than $K$-PC field. Note that every stationary field
over $G$ is  $K$-PC field with $K=G$.
A $K$-PC field is labeled {\em strongly  PC} if $G/K$ is compact, and {\em weakly  PC}
otherwise.

\bigskip
For example if $X$ is a field over $\mbR^2$ (or  $\mbZ^2$) such that  for every $\mfs,\mft \in \mbR^2$ (or $\mbZ^2$),
$\mfK_X(\mfs,\mft) = \mfK_X\big(\mfs+(T_1,0),\mft+(T_1,0)\big)= \mfK_X\big(\mfs+(0,T_2),\mft+(0,T_2)\big)$, $0<T_1, T_2 < \infty$, then $X$ is strongly $K$-PC with $ K = \{(k_1 T_1,k_2,T_2)\!:\! k_1,k_2 \in \mbZ\}$. Since
 $\mfK_X$ is invariant under shifts from $K$ this leads to the existence of unitary operators $U_1$, $U_2$ in $\cH_X$
  such that $U_1 X(\mft) = X(t_1+T_1,t_2)$ and $U_2 X(\mft) = X(t_1,t_2+T_2)$ for every $\mft=(t_1,t_2)$.
If  the field $X$ instead satisfies $\mfK_X(\mfs,\mft) = \mfK_X\big(\mfs+(T_1,T_2),\mft+(T_1,T_2)\big)$, then $X$ is weakly $K$-PC with $ K = \{k(T_1,T_2)\!:\! k \in \mbZ\}$. This leads to a unitary operator $U$ such that $U X(\mft) = X\big(\mft + (T_1,T_2)\big)$,  $\mft=(t_1,t_2)$.

Examples of  PC fields on $\mbZ^2$ can be constructed by a periodic amplitude or time deformation of a stationary field. Suppose $X(\mft)=f(\mft) Y(\mft)$, $\mft=(t_1,t_2)\in \mbZ^2$, where $Y$ is a stationary field and $f$ is a non-random periodic function such that $f(t_1,t_2)=f(t_1+T_1,t_2)=f(t_1,t_2+T_2)$. Then
the field $X$ is strongly $K$-PC with  $ K = \{(k_1 T_1,k_2T_2)\!:\! k_1,k_2 \in \mbZ\}$. If $f$ above instead satisfies $f(t_1,t_2)=f(t_1+T_1,t_2+T_2)$, then $X$ is weakly $K$-PC with $K = \{k(T_1,T_2)\!:\! k \in \mbZ\}$.
 If the function $f$ is two-dimensional integer valued, then the field $X$ defined by $X(\mft)=Y\big(\mft+f(\mft)\big)$ will be weakly PC.
 
\bigskip
Remark that generally $\mfK_X(t+u,s+u) = \mfK_X(t-s+s+u,s+u)$, so a continuous  field $X$ is $K$-PC if and only if   $\mfK_X(t+u,u)$ is a $K$-periodic function of $u$  for every $t\in G$.
If $X$ is $K$-PC field then for all $t,s\in G$  there is a
  unique function $x\mapsto b_X(t,s;x)$ on $G/K$ such that
$$\mfK_X(t+u,s+u)= b_X(t,s; \imath(u)),\quad t,s,u\in G.$$
The canonical map $\imath:G\to G/K$ is continuous and open \cite[Section\,III.1.6]{reiter}, so the function $x\mapsto b_X(s,t;x)$ is continuous. Denote
\begin{equation*} 
\mfB_X(t;x) \defin b_X(t,0;x),\qquad t\in G,\, x\in G/K.
\end{equation*}
Note that $b_X(t,s;\imath(u))=b_X\big(t-s,0;\imath(s+u)\big)= \mfB_X\big(t-s; \imath(s+u)\big)$ for all $t,s,u\in G$.

\bigskip
In this work we need the following notion to proceed to the spectral analysis.
\begin{Definition}\label{def:square-integr}A $K$-PC field $X$ over $G$ is called $G/K$-square integrable if  the function $\mfB_X(0; \cdot)$ is integrable with
respect to the Haar measure on $G/K$.
\end{Definition}

If $X$ is  a $G/K$-square
integrable $K$-PC field, then  from translation-invariance of the Haar measure it follows that for every $t\in G$, $b_X(t,t;\cdot)$ is $\haar_{G/K}$-integrable, and
\begin{eqnarray*}
&&\int_{\xi(G/K)} \|X(t+u)\|_{\cH}^2\, (\haar_{G/K} \circ {\xi}^{-1})(du) =
\int_{G/K} b_X(t,t;x)\, \haar_{G/K}(dx)\\
&&\qquad
= \int_{G/K}\mfB_X(0;x)\, \haar_{G/K}(dx)
=\int_{\xi(G/K)} \|X(u)\|_{\cH}^2\, (\haar_{G/K} \circ {\xi}^{-1})(du) < \infty,
\end{eqnarray*}
where $\xi$ is any  cross-section for $G/K$.
Also note that if $b_X(t,t;\cdot)$ is
$\haar_{G/K}$-integrable for any $t\in G$, then by Cauchy-Schwarz inequality $b_X(t,s;\cdot)$ is $\haar_{G/K}$-integrable for all $t,s\in G$.

\bigskip
When  $X=P$ is a $K$-periodic continuous field, then it is a PC field and we can readilly prove the following equivalence

\vspace{1.5mm}
\centerline{$P$ is $G/K$-square integrable $\Longleftrightarrow$ $P_K\in L^2(G/K;\cH)$,}

\vspace{1.5mm}
\noindent where $P_K$ is the field defined on $G/K$ by $P=P_K\circ\imath$. 

\bigskip
When $X$ 
is a PC process on $\mbR$ with period $T>0$ (i.e. $\mfK_X(t,s) =\mfK_X (t+T,s+T)$ for all $t,s \in \mbR$), it is well known that the SO-spectrum of $X$ can be described as a sequence of complex
measures $\gamma_j$, $j\in \mbZ$,  on $\mbR$ (cf. \cite{HM}). If
 in addition 
$\sum_j \mathrm{Var}(\gamma_j) < \infty$ then the process $X$ is harmonizable and
\begin{equation} \label{1par-spectrum}
\mfK_X(t,s)  = \int\!\!\!\int_{\mbR^2} e^{i(u t-v s)}\,\digamma(du,dv),
\end{equation}
where
$\digamma\!\defin \sum_j \digamma_{\!\!j}$, and $\digamma_{\!\!j}$ is the image of $\gamma_j$  via the mapping $\ell_{j}(u) \defin (u, u-2\pi
j/T)$.
 If $\sum_j \mathrm{Var}(\gamma_j) =\infty$ then $\digamma\!= \sum_j\digamma_{\!\!j}$ can still be viewed as the SO-spectrum of $X$ in the framework of the Schwartz distributions theory (see \cite{mak00}). For more discussion about PC processes please see Example~\ref{ex-TPC} in Section~\ref{sect:Examples}. A corresponding description of  the SO-spectrum is
available for PC sequences ($G=\mbZ$). Let us remark
here that a PC sequence is always harmonizable, but there
are continuous PC processes which are not, see e.g.~\cite{gladyshev61,gladyshev63}.

The above description of the SO-spectrum of a PC process, which originates from Gladyshev's papers \cite{gladyshev61,gladyshev63}, can be easily extended to the case of coordinate-wise strongly periodically correlated fields over $\mbR^n$ or $\mbZ^n$ (see e.g. \cite{Alekseev91,DragJav82,dehayhurd,gaspar04,hurdkf}). The purpose of this work is to describe the SO-spectrum of a $K$-periodically correlated field for any closed subgroup $K$ of an LCA group $G$ and as a particular case when $G = \mbR^m \times \mbZ^n$. We also briefly address the question of  structure of $K$-PC fields.

\section{Covariance Function of a PC Field}
This section contains an extension of Gladyshev's description of the covariance function of one-parameter PC processes (see \cite{gladyshev61,gladyshev63}) to the case of $K$-PC fields. For any $G/K$-square integrable $K$-PC field $X$, define \emph{the spectral  covariance function of the field $X$} (also called cyclic covariance in signal theory, see e.g.~\cite{Gardner06}) by
\begin{equation} \label{alambda}
{a}_{\lambda} (t)\defin \int_{{G/K}} \cha{\lambda,x} \mfB_X(t;x)\, \haar_{G/K}(dx), \quad
\lambda \in \Lambda_K.
\end{equation}

Let $\xi$ be a fixed cross-section for $G/K$. For each $\lambda\in \Lambda_K$ and $t\in G$ let us define an $\cH_X$-valued function $\bZ^{\lambda}(t)$ on $G/K$ by
\begin{equation} \label{Zlambda}
\bZ^\lambda(t) (x)\defin\cha{\lambda,(\imath(t)+x)} X\big(t+\xi(x)\big),
\quad  x \in G/K.
\end{equation}
Notice that $Z^{\lambda}(t)(x)$ depends on the chosen cross-section $\xi$.
From $G/K$-square integrability of $X$ it follows
that for all $\lambda \in \Lambda_K$ and $t\in G$,
$Z^{\lambda}(t)$ is an element of the 
Hilbert space $L^2(G/K;\cH)$.
\begin{Theorem} \label{thm-Corr} Let $X$ be an $\cH$-valued  $G/K$-square integrable $K$-PC field,
and let ${a}_{\lambda}(t)$ and $Z^\lambda(t)(x)$ be as above. Then
the cross-covariance function $ \mfK_Z^{\lambda, \mu}(t,s) \defin
\big(Z^\lambda(t), Z^\mu(s)\big)_\cK$ of the family $\{Z^\lambda(t)\!:\! \lambda
\in \Lambda_K,t\in G\}$
\begin{equation} \label{SCorr}
\mfK_Z^{\lambda,\mu}(t,s) =\cha{\lambda,(t-s)}\, a_{\lambda-\mu}(t-s) =: \mfR^{\lambda,\mu}(t-s).
\end{equation}
If additionally
\begin{description} \label{A}
\item[{\bf [A]}] \qquad \qquad
the function\quad  $\ds G\ni t
\longmapsto a_0(t)$\, is continuous at \ $t=0$,
\end{description}
then $\{Z^\lambda\!:\! \lambda \in \Lambda_K\}$  is a family of jointly stationary fields over $G$  in  $L^2(G/K;\cH_X)$.  
\end{Theorem}

\proof
Let $\xi$ be a fixed cross-section for $G/K$.
Since $X$ is $G/K$-square integrable and $\imath\circ\xi(x) = x$ for any $x\in G/K$, the function  $B(0;\cdot)$ is $\haar_{G/K}$-integrable, so
$Z^{\lambda}(t)\in\cK$ and
\begin{eqnarray*}
\big(Z^\lambda(t), Z^\mu(s)\big)_{\cK}
&=& \int_{G/K} \cha{\lambda,(\imath(t)+x)}
\overline{\cha{\mu,(\imath(s)+x)}}b_X\big(t,s;\imath\circ\xi(x)\big)\, \haar_{G/K} (dx) \\
&=& \int_{G/K} \cha{\lambda ,(\imath(t)+x)}
\overline{\cha{\mu,(\imath(s)+x)}} \mfB_X\big(t-s; \imath(s)+ x\big)\, \haar_{G/K} (dx) \nonumber \\
&=& \cha{\lambda ,(t-s)} \int_{G/K} \cha{(\lambda - \mu), y} \mfB_X(t-s; y)\, \haar_{G/K} (dy)\nonumber \\
&=& \cha{ \lambda,(t-s)}\, a_{\lambda -\mu}(t-s)  \nonumber
\end{eqnarray*}
for all $s,t\in G$.
In view of relation~(\ref{SCorr}), in order to complete the proof  it is enough to
show that for every $\lambda \in \Lambda_K$, the function $G\ni t \to Z^\lambda(t)\in \cK$ is continuous provided  condition {\bf [A]} is satisfied, and this is obvious since by equality~(\ref{SCorr}),
\begin{eqnarray*}
\big\|Z^\lambda(t) -  Z^\lambda(s)\big\|^2_\cK
&= &\mfK_Z^{\lambda,\lambda}(t,t) - \mfK_Z^{\lambda,\lambda}(t,s) -\mfK_Z^{\lambda,\lambda}(s,t) + \mfK_Z^{\lambda,\lambda}(s,s) \\
&=& 2a_0(0) - \cha{\lambda, (t-s)} a_0 (t-s) - \cha{\lambda, (t-s)} a_0 (s-t).\mbox{\qed}
\end{eqnarray*}

\begin{Proposition}
The condition {\bf [A]} in Theorem~\ref{thm-Corr} is satisfied if either

\smallskip
\begin{tabular}{ll}
(i)& $G$ is discrete, or \\
(ii)& $G/K$ is compact, or\\
(iii)&  $X$ is
bounded, and $\mfB_X(0;\cdot)^{1/2}$ is $\haar_{G/K}$-integrable.
\end{tabular}
\end{Proposition}

\proof
Property {\bf [A]} is evident when the group $G$ is discrete.
When  $G/K$ is compact
then $X$ is clearly bounded because $\|X(t)\|_{\cH} = \| X(\xi(x))\|_{\cH}$
where $x = \imath(t) \in G/K$, and $x\mapsto\| X(\xi(x))\|_{\cH}$ is continuous.
Since  $\haar_{G/K}$ is finite,  the continuity of
 the function $\ds t\mapsto a_0(t)=\int_{G/K} \mfB_X(t;x)\, \haar_{G/K}(dx)$
follows therefore from  Lebesgue dominated convergence theorem.

Suppose now that
$\ds \int_{G/K}\mfB_X(0;x)^{1/2} \,\haar_{G/K}(dx) <\infty$.
In this case for all $t, s,x$
\begin{eqnarray*}
&&\big|\mfB_X(t;x)-\mfB_X(s;x)\big|
\,=\,\big|\mfK_X\big(t+\xi(x),\xi(x)\big)-\mfK_X\big(s+\xi(x),\xi(x)\big)\big|\\
&&\qquad\qquad\leq\, \big\|X\big(t+\xi(x)\big)-X\big(s+\xi(x)\big)\big\|_{\cH} \, \big\|X(\xi(x))\big\|_{\cH} \leq 2\sup_t\|X\|_{\cH}
\mfB_X(0,x)^{1/ 2},
\end{eqnarray*}
and
$\ds\lim_{u\to t}\mfB_X(u;x)=\mfB_X(t;x)$. Hence  Lebesgue dominated convergence
theorem applies and we conclude that
$\lim_{s\to t}a_0(s)=a_0(t),$
so  condition {\bf[A]} is satisfied. \qed

\bigskip
Relation~(\ref{SCorr}) in Theorem~\ref{thm-Corr} 
can be also obtained using Gladyshev's technique, that
is by showing non-negative definiteness of
$\ds
\left[\mfR^{\lambda, \mu}(t)\right]_{\lambda,\mu \in \Lambda_K}$, i.e. that
\begin{equation} \label{posdef}
\sum_{j=1}^n \sum_{k=1}^n c_j\overline{c_k}\,\mfR^{\lambda_j,\lambda_k}(t_j - t_k) \geq 0,
\end{equation}
for 
any finite set of complex numbers $\{c_1, \dots, c_n\}$.
Our method, which is an adaptation of the technique used in \cite{mak99},
has the advantage that it gives an explicit construction of an associated stationary family of fields.
We want to point out here that even in the case of PC processes on $\mbR$ with period $T$
 not every matrix function
$\ds \left[\mfR^{m, n}(t)\right]_{m,n \in \mbZ} $
with continuous entries, which is non-negative definite in the sense of~(\ref{posdef})
and such that $ \mfR^{m,n}(t)\,e^{i2\pi  mt/T  }$ depends only on $
m-n$, is associated with a continuous PC process of period $T$
through the relation~(\ref{SCorr}). To achieve the one-to-one
correspondence one has to  consider not necessarily continuous PC
processes (see e.g.~\cite{mak99}).

\bigskip
To complete the analysis of the family of fields $\{Z^{\lambda}\!:\!\lambda\in\Lambda_K\}$ defined by~(\ref{Zlambda}), consider the space $\cH_Z\defin \lin{Z^\lambda(t)\!:\!\lambda \in \Lambda_K, t\in G}$. Clearly $\cH_Z$ is a subspace of $L^2(G/K; \cH_X)$. 
In the case where $G=\mbZ^n\times\mbR^m$, these Hilbert spaces coincide. More precisely
\begin{Proposition}\label{prop:MZ=L2cHX}
Let $X$ be an $\cH$-valued $G/K$-square integrable $K$-PC field, and $Z^{\lambda}$ be defined as above by~(\ref{Zlambda}). Assume that the LCA group $G$ admits a countable dense subset, which is true when $G=\mbZ^n\times\mbR^m$. Then $\cH_Z=L^2(G/K; \cH_X)$.
\end{Proposition}

\proof
We know that $\cH_Z\subset L^2(G/K; \cH_X)$.
 To show the equality, let $f\in L^2(G/K; \cH_X)$ be such that
$\big(Z^{\lambda}(t), f\big)_{\cK} = 0$ for all $\lambda\in \Lambda_K$ and $t\in G$,
i.e.
$$\int_{G/K} \cha{\lambda,(\imath(t)+x)} \big(X\big(t+\xi(x)\big),f(x)\big)_{\cH}\, \haar_{G/K}(dx)= 0, \quad\lambda\in \Lambda_K,\,t\in G.$$
Since $\Lambda_K\sim\what{G/K}$ and  $\cha{\lambda,\imath(t)}\neq 0$,  the scalar product $\big(X(t+\xi(x)),f(x)\big)_{\cH} = 0$, for $\haar_{G/K}$-almost every $x\in G/H$ and  for every $t\in G$~\cite[Theorem 23.11]{HR}.
This implies that for each $t$ there is a negligible Borel subset
$\Xi_t$ of $G/K$ such that $X\big(t+\xi(x)\big)\perp_{\cH} f(x) $ for every
$x\notin \Xi_t$. Note that for every $x$,
$\lin{X\big(t+\xi(x)\big)\!:\! t\in G}= \lin{X(t)\!:\! t\in G} = \cH_X$. If $G$ admits a countable dense
subset  $G^*$, which is true in the case when
$G=\mbZ^m\times\mbR^n$, then from continuity of $X$, it follows
that  also $\lin{X\big(t+\xi(x)\big)\!:\! t\in G^*} = \cH_X$. Therefore
$f(x)\perp_{\cH} \cH_X$ for all $x$ which are not in the negligible set
$\bigcup_{t\in G^*}\Xi_t$.  Hence $f(x)=0$ for $\haar_{G/K}$-almost every $x\in G/K$ and  Proposition~\ref{prop:MZ=L2cHX} is proved.
\qed

\section{SO-spectrum of a PC Field} \label{sec-spectrum}
Let $X$ be a $G/K$-square integrable $K$-PC over $G$ and let $\mfK_X(t,s)$ be its covariance
function. The objective is to describe the domain of the SO-spectrum of $X$, which by definition
(see Definition \ref{def:spectrum of X}) is the domain of the spectrum of the function
$\Phi(t,s) =\mfK_X(t,-s)$.

First we give a description of the domain of the domain of the SO-spectrum of the PC field $X$ in the simplest case where $G/K$ is compact.
\begin{Lemma} \label{lem:domain of PC} Let $X$ be a $G/K$-square integrable $K$-PC over $G$,
and let $\Lambda_K = \imath^*(\what{G/K})\subseteq \what{G}$. Then the domain of the SO-spectrum of
$X$ is the subgroup $L$ of $\what{G}\times \what{G}$ given by
\begin{equation*}
L = \{(\gamma, \gamma - \lambda)\!:\! \lambda \in \Lambda_K, \gamma\in
\what{G}\}.  \end{equation*}
\end{Lemma}
Note that $ L $ can be viewed as the union of {\em hyperplanes}
\begin{equation} \label{FT-K}
L =  \bigcup_{\lambda \in \Lambda_K} L_\lambda\qquad\mbox{where}\quad L_\lambda \defin \{(\gamma,
\gamma - \lambda)\!:\!  \gamma\in \what{G}\}. \end{equation}
\proof
Since $\mfK_X(t+u,s+u)$ is $K$-periodic in $u$,  the function  $\Phi(t,s) =\mfK_X(t,-s)$ is
itself a periodic function on ${G} \times {G}$  with the period $D = \{(k,-k)\!:\! k\in K
\}\subseteq G \times G$. The domain of the SO-spectrum of $X$ is therefore  the dual group of $(G \times
G)/D$ viewed as a subgroup of  $\what{G}\times \what{G}$. 
Note that the subgroup   $D$ is the image of
the subgroup $\{0\} \times K$ through the isomorphism ${\cal D}\!:\! G \times G \ni (t,s) \mapsto  (t+s,-s) \in G \times G$, and this induces an isomorphism
from the quotient group $(G\times G)/D$ onto $(G\times G)/(\{0\}\times K)$.
Futhermore, since $(G \times G)/(\{0\} \times K) = G \times G/K$ and
its dual is $\what{G} \times \Lambda_K$, we deduce that the dual of
$(G \times G)/D$ can be identified with the subgroup $L$ of $\what{G} \times \what{G}$
consisting of the elements of the form $(\chi, \chi- \lambda)$,  $\chi \in \what{G}$, $
\lambda \in \Lambda_K$. 
\qed

\bigskip
We have not used the fact that $\mfK_X$ is a covariance function of a process. 
It turns out that this additional property of  $\mfK_X$  (i.e. the fact that it is nonnegative definite) implies that  the "part of the SO-spectrum" that sits on each $L_\lambda$ is a measure. We want to point out that the set $\Lambda_K$ may be uncountable.

\begin{Theorem} \label{thm-spectrum} Suppose that $X$ is a $G/K$-square integrable $K$-PC field that satisfies the condition {\bf [A]} of Theorem \ref{thm-Corr}.
Then for  every $\lambda \in \Lambda_K$ there is a unique Borel complex measure $\gamma_\lambda$ on $\what{G}$ such
that
\begin{equation} \label{inta}
a_\lambda (t) = \int_{\what{G}}\overline{\cha{\chi, t}}\, \gamma_\lambda(d\chi),
\quad t\in G.
\end{equation}
Furthermore $\sup_\lambda \mathrm{Var}(\gamma_\lambda ) < \infty$ and for each $\lambda\in\Lambda_K$, $\gamma_{\lambda}$ is absolutely continuous with respect to $\gamma_0$.
\end{Theorem}

In this paper by a \emph{representation of $G$ in a Hilbert
space $\cK$} we  mean a weakly continuous  group $\mcU\defin\{U^t\!:\!t\in G\}$  of unitary operators in $\cK$ (see
\cite[Section\,22]{HR}). In this case
there exists a
weakly countably additive orthogonally scattered (w.c.a.o.s) Borel  operator-valued measure
$E$ on $\what{G}$ such that for every Borel set $\Delta$ the
operator $E(\Delta)$ is an orthogonal projection in $\cK$, and for
every $u,v \in \cK$,
$$\big(U^tu,v\big)_{\cK} = \int_{\what{G}} \overline{\cha{\chi,t}}\,\big(E(d\chi)u,v\big)_{\cK}, \quad
t\in G.$$
''\emph{Orthogonally scattered}'' means that $\big(E(\Delta_1)u, E(\Delta_2)v\big)_{\cK}= 0$ for all disjoint $\Delta_1, \Delta_2$ and $u,v\in \cK$. The measure $E$ will be referred to as the \emph{spectral resolution of the unitary operator group $\mcU$}. 

\bigskip
\proof\!\!\!\emph{of Theorem~\ref{thm-spectrum}}.
The joint stationarity of  the fields $\{Z^\lambda\!:\!\lambda \in \Lambda_K\}$ defined by~(\ref{Zlambda}),
implies
that each $Z^\lambda(t) = U^t Z^\lambda(0)$ where $\mcU\defin\{U^t\!:\!t\in G\}$ is the common shift operators group. 
Condition {\bf [A]} guarantees the continuity of the representation $\mcU$ of $G$ in $L^2(G/K;\cH_X)$, and hence
$$ Z^\lambda(t) = \int_{\what{G}} \overline{\cha{\chi, t}}\,  E(d\chi)Z^\lambda(0),$$
where $E$ is the spectral resolution of $\mcU$. 
 Therefore for every $\lambda, \mu \in \Lambda_K$ there is a complex measure $\Gamma^{\lambda, \mu}$ on $\what{G}$ such that
$ \mfR_Z^{\lambda,\mu}(t) = \int_{\what{G}} \overline{\cha{\chi, t}}\,\Gamma^{\lambda,\mu}(d\chi), $
namely, $\Gamma^{\lambda, \mu}(\Delta) = \big(E(\Delta) Z^\lambda(0), Z^\mu(0)\big)_\cK $, where $\cK\defin L^2(G/K;\cH_X)$.  Consequently,  from relations~(\ref{alambda}) and~(\ref{SCorr}) we conclude that
$$ {a}_{\lambda} (t)=\int_{G/K} \cha{\lambda,x} \mfB_X(t;x)\,\haar_{G/K} (dx)  = \mfR_Z^{0, -\lambda}(t) = \int_{\what{G}} \overline{\cha{\chi,t}}\, \Gamma^{0,-\lambda}(d\chi).$$
Then equality~(\ref{inta}) is satisfied with $\gamma_\lambda \defin \Gamma^{0, -\lambda}$.
From Cauchy-Schwarz inequality we have
\begin{eqnarray*}
&&\left|\Gamma^{0, -\lambda}(\Delta) \right| =
\left| \big(E(\Delta) Z^0(0), E(\Delta)Z^{-\lambda}(0)\big)_\cK \right|\\
&&\qquad\qquad\quad\leq \,
\sqrt{\Gamma^{0, 0}(\Delta)}\sqrt{\Gamma^{-\lambda,-\lambda} (\Delta)} = \sqrt{\Gamma^{0,0}(\Delta)} \sqrt{\Gamma^{0, 0}(\Delta -\lambda)},
\end{eqnarray*}
and we deduce the absolute continuity of $\gamma_{\lambda}$ with respect to $\gamma_0$ for any $\lambda \in \Lambda_K$.

Finally note that the total variations of measures $\Gamma^{\lambda, \lambda}$,
$\lambda \in \Lambda_K$, are all equal to $\Gamma^{0,0}(\what{G})$. Indeed, since the measures $\Gamma^{\lambda, \lambda}$,
$\lambda \in \Lambda_K$, are non-negative
$$ \mathrm{Var}\left(\Gamma^{\lambda, \lambda}\right) = \Gamma^{\lambda,\lambda}(\what{G})
=\mfR_Z^{\lambda,\lambda}(0)
=\int_{G/K} \mfB_X(0;y)\,\haar_{G/K}(dy) = \mfR_Z^{0,0}(0) = \Gamma^{0,0}(\what{G}).
$$
Hence all total variations
$\mathrm{Var}\left(\Gamma^{\lambda, \mu}\right)
\leq\sqrt{\Gamma^{\lambda, \lambda}(\what{G})} \sqrt{\Gamma^{\mu,\mu}(\what{G})}$, $\lambda,\mu\in\what{G}$,
are bounded by the same constant, and in consequence all measures $\gamma_\lambda$,
$\lambda \in \Lambda_K$, have uniformly bounded total variations.\qed

\bigskip
Remark that when the field $X=P$ is $K$-periodic and $G/K$-square integrable,   the field $P_K$ is $\haar_{G/K}$-square integrable and  thanks to Parseval equality, the spectral covariance function of the field $P$ can be expressed as
\begin{eqnarray*}
a_{\lambda}^P(t)&=&
\int_{G/K}\cha{\lambda,x}\big(P_K(\imath(t)+x),P_K(x)\big)_{\cH}\,\haar_{G/K}(dx)\\
&=&
\int_{\Lambda_K}\overline{\cha{\chi,\imath(t)}}\big(\what{P_K}(\chi),\what{P_K}(\chi-\lambda)\big)_{\cH}\,\haar_{\Lambda_K}(d\chi)
\end{eqnarray*}
where $\what{P_K}$ is the Fourier Plancherel transform of the field $P_K$.
Then, we deduce that the function $\chi\mapsto \big(\what{P_K}(\chi),\what{P_K}(\chi-\lambda)\big)_{\cH}$ is the density function of the   SO-spectral measure $\gamma_{\lambda}^P$ of the field $P$ with respect to $\haar_{\Lambda_K}$,
\begin{equation}\label{gamma-P}
\gamma_{\lambda}^P(\Delta)=\int_{\Delta\cap\Lambda_K}\big(\what{P_K}(\chi),\what{P_K}(\chi-\lambda)\big)_{\cH}\,\haar_{\Lambda_K}(d\chi)
\end{equation}
for any $\Delta\in\borel(\what{G})$. Particulary
$$\gamma_{0}^P(\Delta)=\int_{\Delta\cap\Lambda_K}\big\|\what{P_K}(\chi)\big\|_{\cH}^2\,\haar_{\Lambda_K}(d\chi).$$
Notice that SO-spectral measure $\gamma_{\lambda}^P$ is concentrated on $\Lambda_K\subset\what{G}$ : $\gamma_{\lambda}^P(\Delta)=\gamma_{\lambda}^P(\Delta\cap\Lambda_K)$ for any $\Delta\in\borel(\what{G})$.
When in addition $\what{P_k}$ is $\haar_{\Lambda_K}$-integrable, then the fields $P_K$ and $P$ are harmonizable with
\begin{eqnarray*}
\mfK_P(t,s)&=&\mfK_{P_K}(\imath(t),\imath(s))\\
&=&\int\!\!\int_{\Lambda_K\times\Lambda_K}\overline{\cha{\lambda,\imath(t)}}\cha{\mu,\imath(s)}\big(\what{P_K}(\lambda),\what{P_K}(\mu)\big)_{\cH}\,\haar_{\Lambda_K}(d\lambda)\haar_{\Lambda_K}(d\mu)\\
&=&\int\!\!\int_{\what{G}\times\what {G}}
\overline{\cha{\chi,t}}\cha{\beta,s}\,\digamma^P(d\chi,d\beta)
\end{eqnarray*}
where $\digamma^{\!P}$ is the measure on $\what{G}\times\what {G}$ concentrated on $\Lambda_K\times\Lambda_K$ defined by
$$\digamma^{\!P}(\Delta)\defin\int\!\!\int_{\Delta\cap(\Lambda_K\times\Lambda_K)}\big(\what{P_K}(\chi),\what{P_K}(\beta)\big)_{\cH}\,\haar_{\Lambda_K}(d\chi)\haar_{\Lambda_K}(d\beta),\quad\Delta\in\borel(\what{G}\times\what{G}).$$
From relation~(\ref{gamma-P}) we find out that the measure $\gamma_{\lambda}^P$ or more precisely its image $\digamma_{\!\!\lambda}^P\defin\gamma_{\lambda}^P\circ \ell_{\lambda}^{-1}$ through the mapping $\ell_\lambda: \what{G} \to \what{G}^2 $ defined by $\ell_\lambda(\chi) = (\chi,
\chi-\lambda)$, $\chi\in \what{G}$, is the restriction of the measure $\digamma^P$ to the hyperplane $L_{\lambda}\defin\{(\chi,\chi-\lambda)\!:\!\chi\in\what{G}\}$.

\bigskip
More generally, when $X$ is a PC field,  the family  $\{\gamma_\lambda\!:\!\lambda\in\Lambda_K\}$  is commonly referred to as  the
\emph{SO-spectral family of} $X$. The measure $\gamma_\lambda$ or more precisely its image $\digamma_{\!\!\lambda}=\gamma_\lambda \circ \ell_\lambda^{-1}$ represents the part of the SO-spectrum of $X$ that sits on the hyperplane $L_\lambda = \{(\chi, \chi-\lambda)\!:\!  \gamma\in \what{G}\}$, see relation~(\ref{FT-K}). 

\bigskip
Next, we give a sufficient condition for a PC field to be harmonizable. 
\begin{Theorem} \label{thm-harmsp} Let $X$ be a $G/K$-square integrable $K$-PC field that satisfies the condition {\bf [A]} of Theorem \ref{thm-Corr}, and let $\{\gamma_\lambda\!:\!\lambda\in\Lambda_K\}$ be
the SO-spectral family of $X$.   Suppose that there is an $\haar_{\Lambda_K}$-integrable non-negative function $\omega$ on $\Lambda_K$ such that for every $\lambda \in \Lambda_K$
\begin{equation} \label{harm-cond}
|\gamma_\lambda(\Delta) | \leq \omega(\lambda)
\quad\mbox{for any Borel }\,\,\Delta \in \borel(\what{G}).
\end{equation}
Then the field $X$ is harmonizable and the SO-spectral measure of $X$ is given by
\begin{equation} \label{F}
\digamma\!(\Delta) = \int_{\Lambda_K}\digamma_{\!\!\lambda}(\Delta)\, \haar_{\Lambda_K}(d\lambda),
\quad\mbox{for any Borel }\,\,\Delta \in \borel(\what{G}\times \what{G}),
\end{equation}
where $\digamma_{\!\!\lambda}\defin\gamma_\lambda \circ \ell_\lambda^{-1}$ and $\ell_\lambda(\chi) \defin (\chi,
\chi-\lambda)$, $\chi\in \what{G}$.
\end{Theorem}

Notice that condition~(\ref{harm-cond}) is satisfied by any  $G/K$-square integrable $K$-periodic field $P$ such that $\what{P_K}$ is $\haar_{\Lambda_K}$-integrable. Here we can take $\omega(\lambda)$ equal to the total variation of the SO-spectral measure 
$\gamma^P_{\lambda}$ of the field $P$
$$\omega(\lambda)=\int_{\Lambda_K}\left|\big(\what{P_K}(\chi),\what{P_K}(\chi-\lambda)\big)_{\cH}\right|\,\haar_{\Lambda_K}(d\chi).$$
The integrability condition on $\what{P_K}:\Lambda_K\to\cH$  is always satisfied when $\Lambda_K$ is compact that is when $G/K$ is discrete, and in particular when $G=\mbZ^n$.

\bigskip
\proof\!\!\emph{of Theorem~\ref{thm-harmsp}}. Let $\{Z^\lambda\!:\!\lambda \in \Lambda_K\}$ be as in Theorem~\ref{thm-Corr}.
From the proof of Theorem~\ref{thm-spectrum} it follows that
$$
\gamma_\lambda(\Delta) =
\Gamma^{0, -\lambda}(\Delta) = \big(E(\Delta) Z^0(0), Z^{\-\lambda}(0) \big)_\cK
$$
Thanks to   definition~(\ref{Zlambda}) and Lebesgue dominated convergence theorem
it follows that the field $\Lambda_K \ni \lambda \mapsto
Z^{\-\lambda}(0)\in \cK$ is continuous, and hence  by assumption~(\ref{harm-cond}), $\lambda \mapsto \gamma_\lambda(\Delta)$ is
integrable over $\Lambda_K$  for every Borel $\Delta$ of $\what{G}$. For all
Borel $D\subseteq \Lambda_K$ and $\Delta\subseteq \what{G}$ let us define
\begin{equation*} 
\tilde{\digamma}(\Delta \times D) \defin \int_D\gamma_\lambda(\Delta)\, \haar_{\Lambda_K}(d\lambda)
= \int_D \big(E(\Delta) Z^0(0), Z^{\-\lambda}(0)\big)_\cK\,\haar_{\Lambda_K}(d\lambda).
\end{equation*}

Condition~(\ref{harm-cond}) and again Lebesgue dominated convergence
theorem entail that the function $\tilde{\digamma}(\Delta \times D)$ is
countably additive in $\Delta$ and $D$ separately. So to show that
the bimeasure $\tilde{\digamma}$ extends to a Borel measure  on $\what{G}\times
\Lambda_K$, it is sufficient to show that its Vitali variation is finite (see~\cite{DS,rao}), that is
$$\sup\left\{\sum_{i=1}^n \sum_{j=1}^n \big|\tilde{\digamma}(\Delta_i \times D_j)\big|\!:\!\Delta_i\cap\Delta_j=\emptyset\,\,\mbox{and}\,\,D_i\cap D_j=\emptyset\,\,\mbox{for}\,\, i\neq j\,\, \mbox{in}\,\, \{1,\dots,n\}\right\}< \infty.$$
Since $\sum_{i=1}^n \big|\gamma_\lambda(\Delta_j)\big| \leq
\mathrm{Var}(\gamma_\lambda) \leq 4 \omega(\lambda)$ and the function $\omega$ is $\haar_{\lambda}$-integrable,
\begin{eqnarray*}
\sum_{i=1}^n \sum_{j=1}^n \big|\tilde{\digamma}(\Delta_i \times D_j)\big|
&\leq &
\sum_{j=1}^n \int_{D_j}\sum_{i=1}^n \big|\gamma_\lambda(\Delta_i)\big|\,\haar_{\Lambda_K}(d\lambda) \\
&\leq & \int_{\bigcup_jD_j} 4 \omega(\lambda)\,\haar_{\Lambda_K}(d\lambda)
\leq \int_{\Lambda_K} 4 \omega(\lambda)\,\haar_{\Lambda_K}(d\lambda) <\infty.
\end{eqnarray*}
Hence $\tilde{\digamma}$ is a measure and in particular 
Fubini and Lebesgue dominated convergence theorems hold for $\tilde{\digamma}$.
Let $\Delta\subseteq \what{G}$ be fixed and let $\varphi(\lambda) \defin
\sum_{j=1}^n b_j 1_{D_j}(\lambda)$ be a simple function on $\Lambda_K$. Then
$$
\int_{\Lambda_K} \varphi(\lambda) \tilde{\digamma}(\Delta, d\lambda)
= \sum_{j=1}^n b_j\int_{D_j} \gamma_\lambda(\Delta) \, \haar_{\Lambda_K}(d\lambda)
=\int_{\Lambda_K} \varphi(\lambda) \gamma_\lambda(\Delta)  \,\haar_{\Lambda_K}(d\lambda).$$
From condition~(\ref{harm-cond}) we deduce that
$ \int_{\Lambda_K} \varphi(\lambda)\, \tilde{\digamma}(\Delta, d\lambda)
=\int_{\Lambda_K} \varphi(\lambda) \gamma_\lambda(\Delta)\, \haar_{\Lambda_K}(d\lambda)$
for any bounded Borel function $\varphi$.
Consequently, for any simple function $\phi$ on $\what{G}$ and bounded $\varphi$ on $\Lambda_K$
\begin{equation} \label{iterated}
\int\!\!\!\int_{\what{G}\times\Lambda_K}\phi(\chi) \varphi(\lambda)\,\tilde{\digamma}(d\chi, d\lambda)
=\int_{\Lambda_K} \varphi(\lambda)\left[\int_{\what{G}}\phi(\chi)\,\gamma_\lambda(d\chi)\right]\haar_{\Lambda_K}(d\lambda).
\end{equation}
If $|\phi|$ is bounded by some finite $c>0$ then by condition~(\ref{harm-cond}), the integral
$\int_{\what{G}} |\phi(\chi)|\,\gamma_\lambda(d\chi)$
is bounded by $4c\,\omega(\lambda)$ which is an $\haar_{\lambda}$-integrable function of $\lambda$. Therefore by  Lebesgue dominated convergence
theorem, relation~(\ref{iterated}) holds for any two bounded measurable functions $\phi$ on $\what{G}$ and  $\varphi$ on $\Lambda_K$.
In particular
\begin{eqnarray} \label{GFT}
\int\!\!\!\int_{\what{G}\times\Lambda_K} \overline{\cha{\chi, t}}\, \overline{\cha{\lambda, x}}\,\tilde{\digamma}(d\chi,d\lambda)
&=& \int_{\Lambda_K} \overline{\cha{\lambda,x}} \left[\int_{\what{G}} \overline{\cha{\chi,t}}\, \gamma_\lambda(d\chi)\,\right]\haar_{\Lambda_K}(d\lambda) \\
&=& \int_{\Lambda_K} \overline{\cha{\lambda,x}}a_\lambda(t)\, \haar_{\Lambda_K}(d\lambda) =  \mfB_X(t; x)= \mfK_X(t+x, x) \nonumber
\end{eqnarray}
for all $t\in G$ and $x\in G/K$.

Let $\ell:  \what{G}\times \Lambda_K  \to \what{G}^{\,2}$ be defined by
$\ell(\chi,\lambda) \defin\ell_{\lambda}(\chi)= (\chi, \chi-\lambda)$, and let
$\digamma = \tilde{\digamma}\! \circ \ell^{-1}$ be the image of the measure
$\tilde{\digamma}\!$ through the mapping $\ell$, that is
$ \digamma\!(\Delta) =\tilde{\digamma}\!\{ (\chi,\lambda)\!:\! (\chi, \chi-\lambda) \in \Delta\}$.
Then $\digamma\!$ is a Borel measure on $\what{G}^{\,2}$ and change of variables formula yields that
\begin{equation}\label{chvar}
\int\!\!\!\int_{\what{G}\times\what{G}} \psi(\chi_1,\chi_2)\, \digamma\!(d\chi_1,d\chi_2) = \int\!\!\!\int_{\what{G}\times\Lambda_K} \psi(\chi,\chi-\lambda)\,\tilde{\digamma}\!(d\chi,d\lambda)
\end{equation}
for any bounded Borel function $\psi:\what{G}\times\what{G}\to\mbC$.
In particular, in view of relation~(\ref{GFT})
\begin{eqnarray*}
\int\!\!\!\int_{\what{G}\times\Lambda_K} \overline{\cha{\chi,t}} \cha{\lambda,s}\,\digamma\!(d\chi, d\lambda)
&=&
\int\!\!\!\int_{\what{G}\times\Lambda_K}\overline{\cha{\chi,t}}\cha{(\chi-\lambda),s}\,
\tilde{\digamma}\!(d\chi, d\lambda) \\
&=& \mfB_X\big(t-s; \imath(s)\big) = \mfK_X(t,s).
\end{eqnarray*}
for all $t,s\in G$ (recall that $\cha{\lambda,s}=\cha{\lambda,\imath(s)}$ for all $\lambda\in\Lambda_K$ and $s\in G$).
Thus the field $X$ is harmonizable and $\digamma\!$ is its SO-spectral measure. Note that
relation~(\ref{iterated}) holds true if the
 product $\phi(\chi) \varphi(\lambda)$ is replaced by any bounded  measurable function $\psi(\chi,\lambda)$ of two variables. Thanks to such upgraded relation~(\ref{iterated}) and to relation~(\ref{chvar}) with $\psi = 1_\Delta$, we get
 \begin{equation*}
 \digamma\!(\Delta)= \int\!\!\!\int_{\what{G}\times\Lambda_K} 1_\Delta(\chi,\chi-\lambda)\,\tilde{\digamma}\!(d\chi,d\lambda)
 =\int_{\Lambda_K} \left[\int_{\what{G}}1_{\Delta}(\chi,\chi-\lambda)\,\gamma_\lambda(d\chi)\right]\haar_{\Lambda_K}(d\lambda)
 \end{equation*}
So, by the definition of $\tilde{\digamma}\!$, we deduce relation~(\ref{F}).
\qed

\bigskip
Note that if $G = \mbZ^n$ then  condition [A]  is satisfied, $\Lambda_K$ is compact, and
condition~(\ref{harm-cond}) holds true with $\omega(\lambda) = \mathrm{Var}(\gamma_\lambda) <\infty$. 
Therefore we generalize the  property of harmonizability of the PC sequences proved in~\cite{gladyshev61}.
\begin{Corollary} \label{cor-harm}
Any $G/K$-square integrable $K$-PC field over $G=\mbZ^n$ is harmonizable.
\end{Corollary}

All the results above simplify significantly if $G/K$ is compact, because then every $K$-PC field over $G$ is $G/K$-square integrable and condition [{\bf A}] in Theorem \ref{thm-Corr} is always satisfied.
\begin{Theorem} \label{Spectrum-compact} Suppose
that $X$ is a $K$-PC field and that $G/K$ is compact. Then for
every $\lambda \in \Lambda_K$ there is a Borel complex measure
$\gamma_\lambda$ on $\what{G}$ such that
$$ a_\lambda (t) =\int_{\what{G}} \overline{\cha{\chi,t}}\, \gamma_\lambda(d\chi),\quad t\in G.$$
Moreover the set $\Lambda_K$ is countable,
$\ds \sum_{\lambda \in \Lambda_K} |a_\lambda (t)|^2  <\infty$
and for every $t\in G$
\begin{equation} \label{series}
\mfB_X(t;x) \stackrel{L^2}{=}
\sum_{\lambda \in \Lambda_K} \overline{\cha{\lambda,x}} a_\lambda (t) \end{equation}
(the series (\ref{series}) converges  in $L^2(G/K)$ with respect to $x$).
Additionally:
\begin{enumerate}
\item[(i)] 
if $\ds\sum_{\lambda\in\Lambda_K} |a_\lambda (t)| <\infty$ for every $t\in G$, then the series~(\ref{series}) converges also pointwise and uniformly with respect to $x\in G/K$, and for all $t,s \in G$ we have
$$\ds\mfK_X(t+s,s) = \sum_{\lambda\in\Lambda_K} \overline{\cha{\lambda,s}} a_\lambda (t);$$
\item[(ii)]
if $\ds\sum_{\lambda\in\Lambda_K} \mathrm{Var}(\gamma_\lambda) < \infty$,
then $X$ is harmonizable, and for all $t,s \in G$ we have
 $$\ds \mfK_X(t,s) = \int\!\!\!\int_{\what{G}\times\what{G}} \overline{\cha{\chi,t}}\, {\cha{\varrho,s}}\,{\digamma}\!(d\chi, d\varrho),
$$
where the SO-spectral measure $\digamma\!$ is given by
$\ds \digamma\!(\Delta) = \sum_{\lambda \in \Lambda_K} \digamma_{\!\!\lambda}(\Delta)$,  $\digamma_{\!\!\lambda} \defin\gamma_\lambda \circ \ell_\lambda^{-1}$ and $\ell_\lambda: \what{G} \to \what{G} \times \what{G} $ is defined as
$\ell_\lambda(\chi) \defin (\chi,\chi-\lambda)$.
\end{enumerate}
\end{Theorem}
\proof Existence of $\gamma_\lambda$ follows from Theorem
\ref{thm-spectrum}. Since for each $t\in G$  the function $x \mapsto
\mfB_X(t;x)$ is bounded, it is in $L^2(G/K)$ and hence its  Fourier
transform $\lambda\mapsto a_\lambda(t)$ is in $L^2(\Lambda_K)$.
Formula~(\ref{series}) is just the inverse formula for
$\mfB_X(t;\cdot)$. Item (i) follows from the uniqueness of the Fourier
transform, while  item (ii) from Theorem \ref{thm-harmsp}. \qed

\section{Structure of PC fields}\label{sect:structure}
When $X$ is a $K$-PC field then for every $k\in K$ the
mapping $V^k: X(t) \mapsto X(t+k)$, $t\in G$,  is well defined and  extends linearly to an isometry from
$\cH_X=\lin{X(t)\!:\! t\in G}$ onto itself. The group $\mcV\defin\{V^k\!:\!k\in K\}$ is a unitary
representation of $K$ in $\cH_X$ and is called the \emph{$K$-shift of $X$}.

\begin{Theorem} \label{PCF-structure}
A continuous field $X$ over $G$ is $K$-PC if and only if  there are a unitary representation $\mcU=\{U^t\!:\!t\in G\}$ of $G$ in $\cH_X$, and a continuous $K$-periodic field $P$ over $G$ with values in $\cH_X$ such that $X(t) = U^t P(t)$, $t\in G$.
\end{Theorem}

\proof
The "if" part is obvious. Prove the other part. Let $X$ be a $K$-PC field, $\mcV=\{V^k\!:\!k\in K\}$ be the $K$-shift of $X$, and $E$ be the spectral resolution of $\mcV$.
Hence $E$ is a w.c.a.o.s. Borel operator-valued measure defined on $\what{K}$.
Since $\what{K}$ is isomorphic to $\what{G}/ \Lambda_K $,
the measure $E$ can be seen as a  measure on $\what{G}/\Lambda_K$,  see \cite[Section\,2.1.2]{rudin}.
Let $\zeta$ be a cross-section for $\what{G}/\Lambda_K$. For every Borel subset
$\Delta$ of $\what{G}$ let us define  $\tilde{E} (\Delta)  \defin E\left( \zeta^{-1}(\Delta)\right)$.
Then  $\tilde{E}$ is a w.c.a.o.s. Borel  operator-valued measure on $\what{G}$ 
whose support is contained in  a measurable set
$\zeta(\what{G}/\Lambda_K)$,  and whose values are orthogonal projections in $\cH_X$.
Since for all $\chi \in \what{K}$ and $k\in K$, $\cha{\chi, k} = \cha{\zeta(\chi),k}$, by change of variable we obtain that
$$ V^k = \int_{\what{G}}  \overline{\cha{\chi,k}}\, \tilde{E}(d\chi), \quad  k \in K.$$
Following Gladyshev's idea (\cite{gladyshev61})
for every $t\in G$ define the operator on $\cH_X$,
$$ U^t \defin \int_{\what{G}} \overline{\cha{\chi,t}}\, \tilde{E}(d\chi), \quad  t\in G.$$
 Clearly $\mcU\defin\{U^t\!:\!t\in G\}$ is a group of unitary operators indexed by
$G$. Moreover for every $v\in \cH_X$,
$$\|(U^t - I)v \|_{\cH}^2 =  \int_{\what{G}} \left| \overline{\cha{\chi,t}} -1 \right|^2\mu_v(dx)$$
where $\mu_v(dx) = \|E(dx)v\|_{\cH}^2$ is a finite non-negative
measure on $\what{G}$. From Lebesgue dominated convergence theorem
we therefore conclude that the unitary operator group 
$\mcU$ is continuous, and hence it is a
unitary representation of $G$ in $\cH_X$.  Note that for $t=k\in K$,
we have  $U^k = V^k$. Define $P(t) \defin U^{-t}X(t)$, $t\in G$. Then
$P$ is continuous and
$$P(t+k) = U^{-t}U^{-k}X(t+k) = U^{-t}V^{-k}X(t+k) = U^{-t}X(t) = P(t), \quad  t\in G,\,k\in K.$$
So $P$ is a continuous $K$-periodic field with values $\cH_X$ and $ X(t) = U^t P(t)$, for every
$t\in G$. \qed

\bigskip
Theorem \ref{PCF-structure} gives a good insight on the origin of the measures $\gamma_{\lambda}$, $\lambda\in\Lambda_K$. Indeed let $X$ be a $G/K$-square integrable $K$-PC field, $\mcU$ and $P$  
be as defined in Theorem \ref{PCF-structure}. Then  $\|X(t)\|_{\cH}=\|P(t)\|_{\cH}$ and $\big(X(t+u),X(u)\big)_{\cH}=\big(U^tP(t+u),P(u)\big)_{\cH}$ for all $t,u\in G$. The PC field $X$ being $G/K$-square integrable, the field $P_K:G/K\to\cH_X$ defined by $P=P_K\circ \imath$ is square integrable  as well as the field $x\mapsto U^tP_K(\imath(t)+x)$ defined on $G/K$,  
for any $t\in G$. Denoting by $\what{P_K}:\Lambda_K\to\cH_X$
the Fourier Plancherel transform of $P_K$,
 the Fourier Plancherel transform of $U^tP_K(\imath(t)+\cdot)$ coincides with the function $$\mu\mapsto\int_{\what{G}}\overline{\cha{\chi+\mu,t}}\,\tilde{E}(d\chi)\what{P_K}(\mu).$$
where $\tilde{E}$ is the   Borel operator-valued measure on $\what{G}$ defined in the proof of Theorem~\ref{PCF-structure}.
Then  thanks to Parseval equality, the spectral covariance function of the PC field $X$ verifies
\begin{eqnarray*}
a_{\lambda}(t)&=&\int_{G/K}\cha{\lambda,x}\big(U^tP_K(\imath(t)+x),P_K(x)\big)_{\cH}\,\haar_{G/K}(dx)\\
&=&
\int_{\Lambda_K}\left(\int_{\what{G}}\overline{\cha{\chi+\mu,t}}\,\tilde{E}(d\chi)\what{P_K}(\mu),\what{P_K}(\mu-\lambda)\right)_{\cH}\,\haar_{\Lambda_K}(d\mu)\\
&=&
\int_{\what{G}}\overline{\cha{\rho,t}}\left(\int_{\Lambda_K}\left(\tilde{E}(d\rho-\mu)\what{P_K}(\mu),\what{P_K}(\mu-\lambda)\right)_{\cH}\,\haar_{\Lambda_K}(d\mu)\right)
\end{eqnarray*}
for all $\lambda\in\Lambda_K$ and $t\in G$.
The SO-spectral measure of the field $X$ is
$$\gamma_{\lambda}(\Delta)=\int_{\Lambda_K}\left(\tilde{E}(\Delta-\mu)\what{P_K}(\mu),\what{P_K}(\mu-\lambda)\right)_{\cH}\,\haar_{\Lambda_K}(d\mu), \quad\Delta\in\borel(\what{G}),\lambda\in\Lambda_K.$$ 
In comparison with expression~(\ref{gamma-P}), we see that the spectral resolution $\tilde{E}$ of the unitary operators group    $\mcU$, in some sense, "spreads" the SO-spectral measure $\gamma_{\lambda}^P$ over $\what{G}$ to form $\gamma_\lambda$.

\bigskip
Theorem~\ref{PCF-structure} also suggests a possibility to decompose a PC field into  stationary components.
If the $K$-PC  field $X$ is  $G/K$-square integrable,  we can therefore \emph{formally}
write
\begin{equation} \label{FTofX}
X(t)  \approx\int_{\Lambda_K} \overline{\cha{\lambda, t}}X^\lambda(t)\,\haar_{\Lambda_K} (d\lambda)
\end{equation}
where  $\{X^\lambda(t)\defin U^t \what{P_K}(\lambda): t\in G\}$, $\lambda\in\Lambda_K $,  is a family of jointly stationary  fields over $G$.
If in addition $\what{P_K}$ is  integrable, then  integral~(\ref{FTofX})  exists, and we have equality for any $t$. 
The integrability condition on $\what{P_K}$   being satisfied if $\Lambda_K$ is compact, that is   in particular when $G=\mbZ^n$, we deduce the following result.

\begin{Corollary} \label{cor-Xrep} Let $X$ be a  $G/K$-square integrable $K$-PC field over $G=\mbZ^n$.
Then there exists a family $\{X^\lambda\!:\! \lambda \in \Lambda_K\}$ of jointly stationary fields over $G$  in  $\cH_X$ such that
\begin{equation*}  
X(t) = \int_{\Lambda_K} e^{i\lambda t'} X^\lambda(t)\, \haar_{\Lambda_K} (d\lambda),
\quad  t\in G,
\end{equation*}
\end{Corollary}
Note that the  pair $(\mcU, P)$ in Theorem \ref{PCF-structure} is
highly non-unique since there are many ways to extend $\mcV=\{V^k\!:\!k\in K\}$ into $\mcU=\{U^t\!:\!t\in G\}$. Consequentlty
 the family $\{X^\lambda\!:\! \lambda \in\Lambda_K\}$ above is likewise not unique.

If $G/K$ is compact, then every $K$-PC field is $G/K$-square
integrable, $\Lambda_K$ is countable, and the integrals above become series.
If $G = \mbZ^n$ and $G/K$ is compact (and hence finite), then  $\Lambda_K$ is finite and Corollary \ref{cor-Xrep}   yields the following $\mbZ^n$ version of
Gladyshev's representation of PC sequences included in  \cite{gladyshev61}.

\begin{Corollary}
Suppose that $X$ is a $K$-PC field over $\mbZ^n$ and that $\mbZ^n/K$ is compact.
Then $\Lambda_K$ is finite and there is a finite family $\{X^\lambda\!:\! \lambda \in \Lambda_K\}$ of jointly stationary fields over $G$  in  $\cH_X$ such that for every $t\in G$
$$X(t) = \sum_{\lambda \in \Lambda_K} e^{i\lambda t'}X^\lambda(t).$$
\end{Corollary}
If $\Lambda_K$ is not compact then even in the case of a periodic
function, its Fourier transform does not have to converge everywhere.

\section{Examples}\label{sect:Examples}

In this section  $G=\mbZ^n\times \mbR^m$, the Haar measure on  $\mbZ^n$ is the counting measure, the Haar measure on  $\mbR^m$ is  $dt/(\sqrt{2\pi})^{m}$ where $dt$ is the Lebesgue measure on  $\mbR^m$, $\what{\mbZ^n}$ will be
identified with  $[0,2\pi)^n$ with addition \emph{mod} $2\pi$, $\what{\mbR^m}$ will be identified with $\mbR^m$, the Haar measures on  $[0,2\pi)^n$ is  $dt/(2\pi)^{n}$, $\what{G} = [0,2\pi)^n \times \mbR^m$,  elements of $G$ and $\what{G}$ are row vectors, and the value of a character $\chi \in \what{G}$ at $t\in G$ is
$\cha{\chi,t} =  e^{-i \lambda  t'}$, where $t'$ is transpose of $t$.  Remembering previous sections, in order to describe the domain of the  SO-spectrum of a $K$-PC field $X$ over $G$ the only task is to identify $G/K$ and $\what{G/K}$ as concrete subsets $Q$ and $\Lambda_K$ of $G= \mbZ^n\times \mbR^m$ and $\what{G} = [0,2\pi)^n\times \mbR^m$, respectively, in the way that the value of character $\lambda \in \Lambda_K$ at $t\in Q$ is still $\cha{\lambda,t} =  e^{-i \lambda  t'}$. This identification, which is obvious when $X$ is coordinate-wise PC, may be less trivial in the case of more complex $K$. It may be helpful, and is worth, to note that any closed nontrivial subgroup $K$ of $G=\mbZ^n\times \mbR^m$ is isomorphic to $\mbZ^k\times \mbR^l$ for some $k, l\in\mbN$ such that $l\leq m$ and  $1\leq k+l\leq n+m$.  This isomorphism, which at least in the case of $G=\mbZ^n$ or $G= \mbR^m$ can be found by selecting a proper basis for $G$ (see \cite[Theorem 9.11 and A.26]{HR}),  provides a description and a parametrization of the sets $Q$ and $\Lambda_K$. As before $\big[a\big]_{b}$ will denote the remained in integer division of $a$ by $b$, $b>0$.

First we briefly revisit a one-parameter case and its slight extension.

\begin{Example} \label{ex-TPC}  {\rm
Suppose that $X$ is a PC process with period $T>0$. Then $K = \big\{kT\!:\! k\in \mbZ\big\}$, $\mbR/K = [0,T)$ with addition modulo $T$ and $\Lambda_K = \big\{ \frac{2\pi k}{T}\!:\! k \in \mbZ\big\}$. Moreover for every $\lambda = \frac{2\pi k}{T}\in \Lambda_K$, $ a_\lambda (t) \defin a_k(t)  = \frac{1}{T} \int_0^{T} e^{-i s \frac{2\pi k}{T}} \mfK_X(t+s,s)\, ds$ and there is a measure $\gamma_k$ on $\mbR$ such that $ a_k(t) = \int_{\mbR} e^{itu}\, \gamma_k(du)$ (Theorem \ref{thm-spectrum}). The domain of the SO-spectrum of $X$ is $L = \bigcup_{k\in \mbZ} L_k$, where $L_k \defin \big\{\big(u, u-\frac{2\pi k}{T}\big)\!:\! u \in \mbR \big\}$. The part of the SO-spectrum that sits on $L_k$ is a measure $\digamma_{\!\!k}$ defined as $\digamma_{\!\!k}=\gamma_k\circ\ell_{k}^{-1}$ where $\ell_{k}: \mbR \to \mbR^2$,   $\ell_{k}(u) = \big(u, u-\frac{2\pi k}{T}\big)$. If $\sum_k \mathrm{Var}(\gamma_k) < \infty$ then the process $X$ is harmonizable and
$\digamma\! \defin \sum_k \digamma_{\!\!k}$ is a measure on $\mbR^2$ which satisfies relation~(\ref{1par-spectrum}).

If $X$ is  stationary  then $K = \mbR$, $\mbR/K = \{0\}$,  $\Lambda_K = \{0\}$,  $a_0 (t) = \int_{\{0\}} \mfK_X(t+s,s)\, \delta_0(ds) = \mfK_X(t,0)$. By Theorem \ref{thm-spectrum} there is a measure $\gamma_0$ on $\mbR$ such that $a_0 (t) = \int_{\mbR} e^{itu}\, \gamma_0(du)$.  Consequently, the SO-spectrum of $X$ is the measure $\digamma\!\! = \digamma_{\!\!0} = \gamma_0\circ \ell_{0}^{-1}$, which sits on the diagonal $L_0 = \big\{(u, u)\!:\! u \in \mbR \big\}$.

To see the need for the square integrability assumption, let us add one parameter to the above process; that is, let us consider a field $X=\{X(s,t)\!:\! (s,t)\in\mbR^2\}$ such
that $\mfK_X\big((s,t), (u,v)\big)  = \mfK_X\big((s+T,t), (u+T,v)\big)$, $s,t,u,v \in \mbR$  ($T>0$ is fixed). Then $K = \big\{(kT,0)\!:\! k\in \mbZ\big\}$, $\mbR^2/K = [0,T)\times \mbR$ with addition modulo $T$ on the first coordinate, $\Lambda_K = \big\{ \big(\frac{2\pi k}{T},x\big)\!:\! k \in \mbZ, x\in \mbR\big\}$ and 
$$a_\lambda (s,t) \defin a_{k,x}(s,t)  = \frac{1}{T\sqrt{2\pi}} \int_0^{T}\!\! \int_{\mbR} e^{-i \big(u\frac{2\pi  k}{T} + vx\big) } \mfK_X\big((s+u,t+v), (u,v)\big)\,  du dv$$
for $\lambda=\big(\frac{2\pi k}{T},x\big)\in \Lambda_K$. 
The square integrability assumption $\int_0^{T} \left[ \int_{\mbR} \| X(s,t)\|_{\cH}^2 dt \right]\,ds < \infty$  assures that the above integral exists. If it does then, in view of Theorem \ref{thm-spectrum}, for every $k\in \mbZ$ and $x\in \mbR$  there exists a measure $\gamma_{k,x}$ on $\mbR^2$ such that $a_{k,x}(s,t)  = \int_{\mbR} e^{i (su + tv) }\, \gamma_{k,x}(du,dv)$.  The domain of the SO-spectrum of $X$ is $L = \bigcup_{k\in \mbZ} \bigcup_{x\in \mbR} L_{k,x}$, where $L_{k,x}$ is a two-dimensional plane in $\mbR^4$, 
$L_{k,x} \defin
\big\{\big(u, v, u-\frac{2\pi k}{T}, v-x\big)\!:\! u,v \in \mbR \big\}$. The "part" of the SO-spectrum that sits on $L_{k,x}$ is a measure $\digamma_{\!\!k,x}$ defined as $\digamma_{\!\!k,x}\defin \gamma_{k,x}\circ\ell_{k,x}^{-1}$, where $\ell_{k,x}: \mbR^2 \to \mbR^4$ is defined by  $\ell_{k,x}(u,v) \defin \big(u, v, u-\frac{2\pi k}{T}, v-x\big)$.  If  $\mathrm{Var}(\digamma_{\!\!k,x}) \leq \omega(k,x)$ and $\sum_k \int_{\mbR} \omega(k,x)\, dx < \infty$, then $X$ is harmonizable and the SO-spectral measure of $X$ is $\digamma\! = \frac{1}{\sqrt{2\pi}} \sum_k \int_{\mbR} \digamma_{\!\!k,x}\, dx$ (see Theorem \ref{thm-harmsp}). Note that $L$ above is, in fact, the union of countably many three-dimensional hyperplanes $D_k$ in $\mbR^4$, $D_k \defin \bigcup_{x\in \mbR} L_{k,x} =\big\{(u, v, u-\frac{2\pi}{T}, v-x)\!:\! u,v,x \in \mbR\big\}$, which are parallel to the "diagonal" $D_0$.

If the field $X=\{X(s,t)\!:\!(s,t)\in\mbR^2\}$ is stationary in $s$, then $\Lambda_K = \big\{ (0,x)\!:\! x\in \mbR\big\}$,  the condition of the square integrability of $X$ means that $ \int_{\mbR} \|X(0,t)\|^2_{\cH}\, dt < \infty$ and, if the latter is satisfied,  the SO-spectrum of $X$ sits on the three-dimensional hyperplane in $\mbR^4$,
$L \defin D_0 = \big\{(u, v+x, u, x)\!:\! u,v,x \in \mbR \big\}$.  \qed
} \end{Example}

Next example contains a complete analysis of the SO-spectrum of a weakly PC field.

\begin{Example} \label{ex-wk} {\rm
Let  $T$ and $S$ be two non-zero integers. Suppose that 	the field $X$ on $\mbZ^2$ is {weakly PC}
with period $(T,S)$   (\cite{hurdkf}), that is
$$\mfK_X\big((m,n),(u,v)\big) = \mfK_X\big((m+T,n+S),(u+T,v+S)\big), \quad\mbox{for all}\,\, n,m,u,v\in\mbZ $$
Here  $K=\big\{k(T,S)\!:\! k\in \mbZ\big\}$. We assume that at least one of $T$ or $S$ is positive.
Let $d\defin \gcd(T,S)$   be the greatest common positive integer divisor of $T$ and $S$, so that $(T,S) = d\times (T_1,S_1)$ and $\gcd(T_1,S_1) = 1$. From Bezout's lemma there are integers $q,p$ such that $T_1 q - S_1 p$ = 1.  Let $\phi$ be a mapping of $\mbZ^2$ onto itself, given by $\phi(m,n) = (m,n)\Phi'$, where
 $\Phi = \left(\begin{array}{cc}
T_1 & p \\ S_1 & q \\
\end{array}  \right)$, and $\Phi'$ stands for the transpose matrix of the matrix $\Phi$.
Because $\det \Phi = 1$, $\phi$ is an isomorphism. Since $\phi(dk,0) = (kT,kS)$ for $k \in \mbZ$, we have $K = \phi(d\mbZ \times \{0\})$ and we identify $G/K$ to $$ Q \defin \phi\big(\{0,\dots, d-1\} \times \mbZ\big) = \big\{(kT_1+lp, kS_1+lq)\!:\! k=0,\dots,d-1, l\in
\mbZ\big\}.$$
The dual mapping $\psi(s,t) = \big[(s,t)  \Phi^{-1}\big]_{2\pi} = \big(\big[qs - S_1t\big]_{2\pi}, \big[- ps  + T_1 t\big]_{2\pi}\big)$, $s, t \in [0,2\pi)$, maps $\big\{\frac{2\pi k}{d}\!:\! k=0,\dots,d-1\big\} \times [0,2\pi)$, which is the dual of $\{0,\dots, d-1\} \times \mbZ$,  onto  the dual $\Lambda_K$ of $Q$. 
The construction that we use produces a  convenient parametrization of  $\Lambda_K$ as the union of $d$ lines : $\Lambda_K = \bigcup_{k=0}^{d-1}\Lambda_k$ 
where
$$ \Lambda_k \defin\left\{\left(\left[\frac{2\pi kq}{d} - S_1t\right]_{2\pi},\left[\frac{- 2\pi kp}{d}  + T_1 t\right]_{2\pi}\right)\!:\! t \in [0,2\pi)\right\}.$$
 Note that the value of a character $(u,v)=\psi(s,t) \in \Lambda_K$ at  $(m,n) = \phi(k,l) \in Q$  is equal to $e^{-i (mu+nv)} = e^{-i (s,t)  \Phi^{-1} \Phi (k,l)'} = e^{-i (ks + lt)}$,  as required. Assume that $X$ is $G/K$-square integrable, for example that  $ \sum_{n = -\infty}^{\infty} \|X(m,n)\|^2_{\cH} < \infty$, for any $m=1,\dots,T-1$. Then from the previous discussion and the results of Section \ref{sec-spectrum}
we deduce the following properties. 
\begin{enumerate} 
\item[(i)]
If $(u,v) \in \Lambda_K$, then $u= \big[\frac{2\pi kq}{d} - S_1t\big]_{2\pi}$, $v= \big[\frac{- 2\pi kp}{d}  + T_1 t\big]_{2\pi}$ for some unique $t\in[0,2\pi)$ and unique $k=0,\dots,d-1$. Hence the spectral covariance $a_{(u,v)} (m,n) =: a_{k,t} (m,n)$ of $X$ at $(m,n) \in \mbZ^2$ is equal to
 \begin{equation*} 
 \hspace{-2mm}
 a_{k,t} (m,n) =   \frac{1}{d} \sum_{j=0}^{d-1} \sum_{l = -\infty}^{\infty} e^{-i\big(j\frac{2\pi k q }{d} + lt\big)} \mfK_X\big((m+jT_1+lp, n+jS_1+lq),(jT_1+lp, jS_1+lq)\big).
 \end{equation*}
\item[(ii)] 
For each $k=0,\dots,d-1$ and $t\in [0,2\pi)$ there exists a measure $\gamma_{k,t}$ on $[0,\pi)^2$ such that
$$ a_{k,t} (m,n) = \int_0^{2\pi} \int_0^{2\pi} e^{i(mx+ny)}\, \gamma_{k,t} (dx, dy). $$
\item[(iii)] 
The SO-spectrum of $X$ sits on the set  $ L = \bigcup_{k=1}^{d-1}  \bigcup_{t\in [0,2\pi)} L_{k,t}$, where $ L_{k,t} $ is a two-dimensional plane in $[0,2\pi)^4$
$$ L_{k,t} \defin\left \{ \left(x,y,\left[x-\frac{2\pi kq}{d} + S_1t\right]_{2\pi},\left[y+\frac{ 2\pi kp}{d} - T_1t\right]_{2\pi}\right)\!:\!
 x, y \in [0,2\pi)\right \}. $$
Note that $D_k =  \bigcup_{t\in [0,2\pi)} L_{k,t}$ is a three-dimensional hyperplane in $[0,2\pi)^4$, so $L$ is, in fact, the union of $d$ three-dimensional hyperplanes.
If $d=\gcd(T,S)=1$, then $L = D_0 = \left \{ \left(x,y,\big[x+St\big]_{2\pi},\big[y-Tt\big]_{2\pi}\right) \!:\!  x, y, t \in [0,2\pi) \right\}.$ In this case the field  $X$ is a rotation of the field  $Y$ defined by $Y(m,n)\defin X\big((m,n) \Phi'\big)$,  $m,n\in\mbZ$, which is stationary in $m$. Indeed $X(m,n)= Y\big((m,n)(\Phi')^{-1}\big)$.
\item[(iv)] 
If there is an integrable function $\omega:[0,2\pi)\to[0,\infty)$ such that  $ Var (\gamma_{k,t}) \leq \omega(t)$ for all $k$ and $t$, then $X$ is harmonizable and
    $$  \mfK_X\big((m,n),(j,r)\big) =  \int_0^{2\pi} \int_0^{2\pi} \int_0^{2\pi} \int_0^{2\pi} e^{i (mu + nv - jx - ry)}\,\digamma\! (du,dv,dx,dy). $$
   The SO-spectral measure  $\digamma\!$ of $X$ is given by $ \digamma\! (\Delta) = \sum_{k=0}^{d-1}  \int_0^{2\pi} \digamma_{\!\!k,t} (\Delta)\, dt $,
 where $\digamma_{\!\!k,t}$ is the complex measure on $[0,\pi)^4$ whose support is contained in the plane $L_{k,t}$ and defined by
$$\hspace{-9mm}\digamma_{\!\!k,t}(\Delta)\!\defin\gamma_{k,t} \left\{ (x,y)\in [0,2\pi)^2\!:\! \left(x,y,\left[x-\frac{2\pi kq}{d} + S_1t\right]_{2\pi},\left[y+\frac{2\pi kp}{d} - T_1t\right]_{2\pi}\right)\in \Delta \right\}. $$
\end{enumerate}
Figure 1 is the graph of the set $\Lambda_K$ defined previously  
 in the case when  $T=12$ and $S = 9$ ($d=3$, $p=q=1$). Then $\Lambda_K=\Lambda_0\cup\Lambda_1\cup\Lambda_2$ consists of three lines, which are shown with different width pattern. If now from each point on the graph we draw the rectangle $[0,2\pi)\times [0,2\pi)$ then the resulting three-dimensional body in $[0,2\pi)^4$ is the domain of the SO-spectrum of the field $X$. \qed
\begin{center}\label{fig1}
\includegraphics[width=2.5in,height=2.5in]{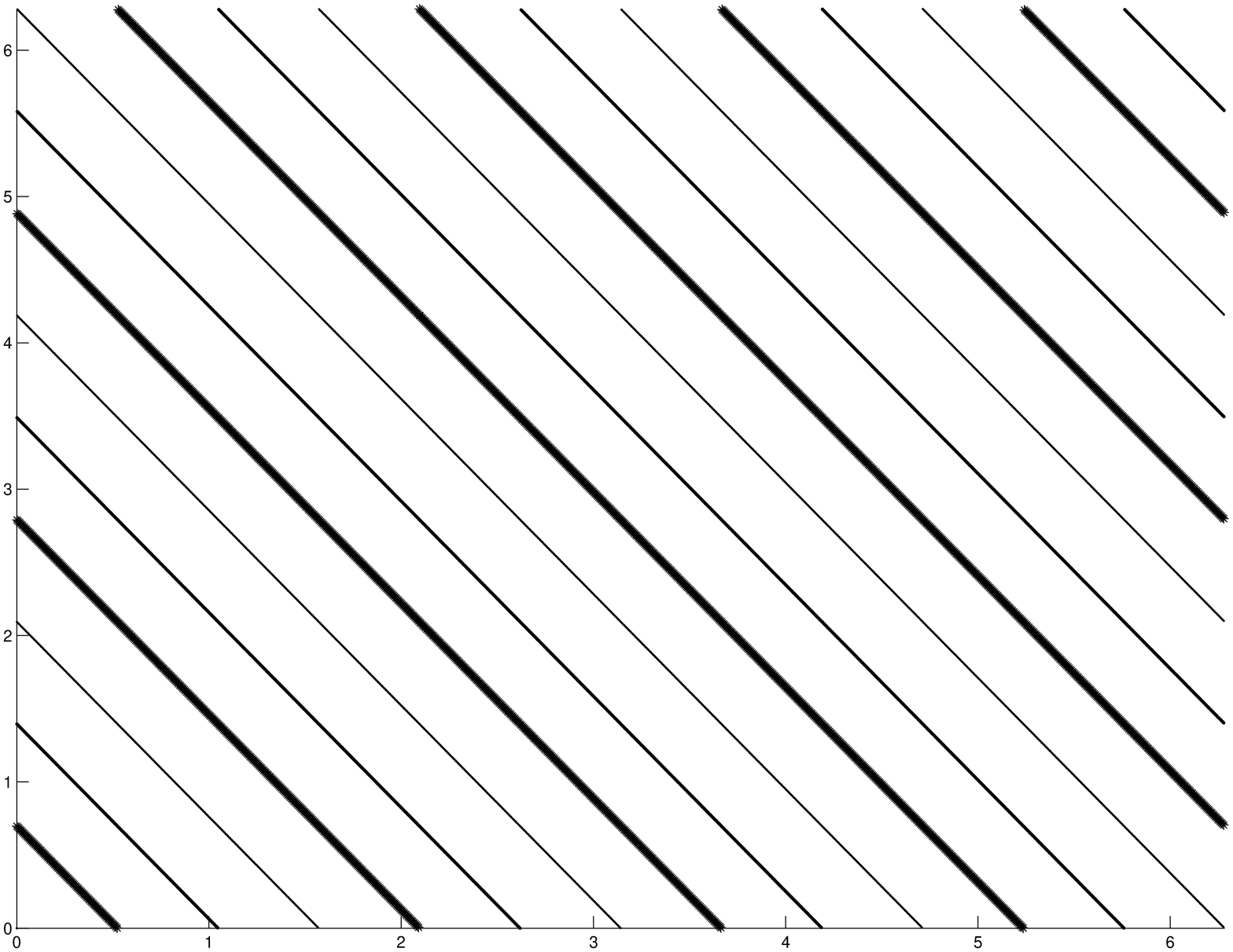} \\
Fig 1: Graph of $\Lambda_K$
\end{center}
}\end{Example}

The last example  is  a particular example of strongly PC fields over $\mbR \times \mbZ^2$.  It combines a mixture of continuous and discrete structures.
\begin{Example}{\rm 
Suppose that $X$ is a field over $G \defin \mbR\times\mbZ^2$ such that
$X(t,m,n) = X(t+4,m,n) = X(t,m+1,n+3) = X(t,m+2,n)$ for all $t\in \mbR$,  $m,n\in \mbZ.$
Then $K = \{k(4,0,0)+j(0,1,3) + l(0,2,0)\!:\! k,l,j \in \mbZ\}$. In order to describe $G/K$ and $\Lambda_K$ we consider a change of basis of $G=\mbR\times\mbZ^2$ defined by  the mapping
$$ \phi(t,m,n) \defin (t,m,n)\Phi',\qquad
\mbox{where}\quad\Phi = \left(
\begin{array}{ccc}
1& 0 & 0 \\
0 &1 & 0\\
0 &3 & 1 \\
\end{array}\right). $$
Then the mapping $\phi$ is an isomorphism of $G$ onto itself and $K = \phi(P)$, where $P= \{(4k, j, 6l)\!:\! k,l,j \in \mbZ\} = 4\mbZ \times \mbZ \times 6\mbZ$. To see this note that $2(0,1,3)  - (0,2,0)= (0,0,6)$, so that $K$ is generated by the 3-tuples $(4,0,0)$, $(0,1,3)$  and  $(0,0,6)$, which are  respectively equal to   $\phi(4,0,0)$, $\phi(0,1,0)$, and $\phi(0,0,6)$. The quotient $G/P = [0,4)\times \{0\} \times \{0,\dots,5\}$, so we take
$Q\defin \phi(G/P) = \left\{ (s,0,l)\!:\! s\in [0,4), l=0,\dots 5\right\}$. The dual of $G/P$ is $\Lambda_P = \frac{2\pi}{T} \mbZ \times \{0\} \times \big\{\frac{\pi r}{3}\!:\! r=0,\dots,5\big\} $ and hence the dual of $G/K$ can be represented as $\Lambda_K = \psi(\Lambda_P)$, where $\psi$ is the isomorphism of $\what{G} = \mbR \times [0,2\pi)^2$ onto itself defined by $\psi(t,u,v) \defin (t,u,v)\Phi^{-1} = \big(t, \big[u-3v\big]_{2\pi}, v\big)$. Therefore
$$ \Lambda_K = \left\{\left(\frac{2\pi k}{T}, \big[-\pi r\big]_{2\pi},  \frac{\pi r}{3}\right)\!:\! k\in \mbZ, r=0,\dots,5 \right\},$$
is countable.
Note that  $\big[-\pi r\big]_{2\pi}$ is either $\pi$ (if $r$ is odd) or 0. For each  $\lambda_{k,r} = \big(\frac{2\pi k}{T}, \big[-\pi r\big]_{2\pi},  \frac{\pi r}{3}\big) \in \Lambda_K$, the corresponding spectral covariance is given by
$$ a_{k,r} (t,m,n) = \frac{1}{24}\sum_{l=0}^5\int_0^4  e^{-i(\frac{2\pi k s}{T}  + \frac{\pi r l}{3})} \mfK_X\big((t+s, m, n+ l),(s, 0, l)\big)\, ds, $$
and for each ${k,r}$ there exists a measure $\gamma_{k,r}$ on $\mbR \times [0,2\pi)^2$ such that $$a_{k,r} (t,m,n) = \int_{\mbR} \int_0^{2\pi} \int_0^{2\pi}e^{i(ts+mu+nv)}_, \gamma_{k,r}(ds, du,dv).$$ The SO-spectrum of the field $X$ sits on the union of countably many hyperplanes
$$L_{k,r} \defin\left\{\left(s,u,v,s-\frac{2\pi k}{T},\left[u+\pi r\right]_{2\pi}, \left[v-\frac{\pi r}{3}\right]_{2\pi}\right)\!:\! s\in \mbR, u,v \in [0,2\pi)\right\},$$
 $k\in \mbZ$, $r=0,\dots,5$, of $\mbR \times [0,2\pi)^2\times\mbR \times [0,2\pi)^2 $. If the sum of total variations of measures $\gamma_{k,r}$ is finite, then $X$ is harmonizable and its SO-spectral measure $\digamma\! = \sum_{k=-\infty}^{\infty} \sum_{r=0}^{5} \digamma_{\!\!k,r}$, where $\digamma_{\!\!k,r} = \gamma_{k,r}\circ \ell_{k,r}^{-1}$ and $\ell_{k,r}(s,u,v) \defin \big(s,u,v,s-\frac{2\pi k}{T},\big[u+\pi r\big]_{2\pi}, \big[v-\frac{\pi r}{3}\big]_{2\pi}\big)$.

Note that if  we define $Y(t,n,m) \defin X\big((t,n,m) \Phi'\big)$,  then $Y$ is PC in $t$ with period $T=4$, stationary in $n$, and  PC in $m$ with period $M=6$. Since  $X(t,n,m) = Y\big((t,n,m) (\Phi')^{-1}\big)$, one can therefore say that the field  $X$ is periodically correlated in direction $(1,0,0)$ with period $T=4$, stationary in direction of $(0,1,3)$ and periodically correlated in direction $(0,0,1)$ with period $M=6$. \hfill \qed
}
\end{Example}

\paragraph{{\bf}ACKNOWLEDGEMENTS} The paper was partially written during the author's stay at Universit\'{e} Rennes 2, Rennes, France, in June 2011.


\bibliographystyle{plain}

\end{document}